\documentclass{article}
\usepackage{chicago}
\renewcommand{\cite}{\shortciteN}
\usepackage{pdfsync}
\usepackage{color}
\usepackage[intlimits]{amsmath}
\usepackage{amsopn}
\usepackage{amsfonts}
\usepackage{amsthm}
\usepackage{fancybox}
\usepackage{graphicx}
\usepackage{bm} 
\usepackage{dsfont} 

\newcommand{\stsets}[1]{\mathbb{#1}}
\newcommand{\A}{\stsets{A}}

\newcommand{\M}{\stsets{M}}
\newcommand{\N}{\stsets{N}}

\newcommand{\R}{\stsets{R}}
\newcommand{\Z}{\stsets{Z}}

\newenvironment{Abs}{\vspace*{0.5cm}\begin{center} {\bf
    Abstract}\end{center}\vspace*{0cm}\nopagebreak 
   \small  } {\vspace{0.5cm}} 

\newenvironment{Key}{\vspace*{0.5cm}\noindent {\bf
    Keywords:}
   \sc } {\vspace{0.5cm}} 

\newenvironment{AMS}{\vspace*{0.5cm}\noindent {\bf
    AMS 1991 Subject Classification.} } {\vspace{0.5cm}}

\theoremstyle{definition}
\newtheorem{Dfn}{Definition}

\theoremstyle{remark}

\newtheoremstyle{mytheorem}{0.5cm}{0.2cm}{\slshape}{ }{\bfseries}{.}{ }{}
\theoremstyle{mytheorem}
\newtheorem{Th}{Theorem}

\DeclareMathOperator{\supp}{supp}

\newcommand{\Unif}{\ensuremath{\mathsf{Unif}}}

\renewcommand{\P}{\mathbf{P}}

\DeclareMathOperator{\E}{{\bf E}}

\newcommand{\cond}{\mid}
\newcommand{\ssp}{\hspace{2pt}}

\newcommand{\mydot}{{\raisebox{.3ex}{$\scriptscriptstyle{\,\bullet\,}$}}}

\newcommand{\mytilde}{{\!\raisebox{-0.9ex}{$\tilde{\ }$}}}

\renewcommand{\epsilon}{\varepsilon}
\renewcommand{\phi}{\varphi}

\newcommand{\thru}{,\dotsc,}
\newcommand{\seg}{{see, e.~g., \nolinebreak}}

\newcommand{\ti}{\to\infty}

\newlength{\querylen}
\newcommand{\ie}{\textit{i.\ssp e.\ }}
\newcommand{\eg}{\textit{e.\ssp g.\ }}

\newcommand{\cf}{\textit{cf.\ }}

\newcommand{\bbb}{\mathcal{B}}

\begin{document}
\title{Optimal design of dilution experiments under volume
  constraints} 

\author{Maryam Zolghadr\thanks{University of Gothenburg, Department of
    Mathematical Sciences, 412 96 Gothenburg, Sweden. Email:
    zolghadr@student.chalmers.se} \and Sergei Zuyev\thanks{Chalmers
    University of Technology, Department of Mathematical Sciences, 412
    96 Gothenburg, Sweden. Email: sergei.zuyev@chalmers.se}}

\maketitle

\begin{Abs}

  The paper develops methods to construct a one-stage optimal design
  of dilution experiments under the total available volume constraint
  typical for bio-medical applications. We consider various design
  criteria based on the Fisher information both is Bayesian and
  non-Bayasian settings and show that the optimal design is typically
  one-atomic meaning that all the dilutions should be of the same
  size. The main tool is variational analysis of functions of a
  measure and the corresponding steepest descent type numerical
  methods. Our approach is generic in the sense that it allows for
  inclusion of additional constraints and cost components, like the
  cost of materials and of the experiment itself.
\end{Abs}

\begin{Key}
  Dilution experiment, optimal design of experiments, Fisher
  information criteria, gradient methods, design measure, variational
  analysis on measures, stem cells counting.
\end{Key}


\section{Introduction}
\label{sec:Introduction}
The paper studies a wide class of statistical experiments with the
total volume constraint arising, in particular, in stem cells
research, a very active area of experimental biology. Stem cells are
the cells produced during early stages of embryonic development and
they are having capacity to turn into various types of tissue
cells. This potentially opens new ways to cure many diseases and
explains the huge importance of the stem cells research, \seg
\cite{Mayhall:2004} and the references therein. We aim to characterise
an optimal design of dilution-type experiments under the contraint of
the available solution which is typical for experiments involving
counting stem cells. Specifically, this study originates in studies of
\textit{Hematopoietic} or \textit{blood stem cells (HSCs)} which are
the stem cells giving rise to all red and white blood cells and
platelets. HSCs develop in a mammal's embryo in early days from
cells. There is not much known on the biological mechanism which
triggers a cell to develop into an HSC and there is no direct way so
far to observe such cells, called \textit{pre-HSCs}, which soon to
become HSCs. Pre-HSCs are mostly produced in aorta-gonad-mesonephros
(AGM) region of the embryonic mesoderm and also in the yolk sac, then
colonise the liver. A challenging problem is how to detect the
pre-HSCs and how many of them are present at a different embryo ages.

In order to estimate the number of pre-HSCs present in a given region,
experiments on laboratory mice have been conducted which use the
following signature property of stem cells. A mature HSC is capable to
cure a mouse which receives a controlled dose of radiation if injected
in its blood. A mouse recovers (\emph{repopulates}), if and only if,
it has received at least one HSC in the injected dose\footnote{This
  assumption is rather questionable, so the biologists cautiously
  speak of one \emph{repopulation unit} for this unknown minimal
  number of HSCs sufficient to cure an irradiated mouse. In this study
  we basically loosely speak of \emph{one HSC} as of \emph{one
    repopulation unit}.}. Thus the number
of pre-HSCs can be estimated by the so-called \textit{limiting
  dilution method}: controlled dozes of a substrate containing samples from the
AGM are injected into irradiated mice and then the number of repopulated
mice infers on the number of HSCs which developed from the initially
present pre-HSCs. 

The dilution method has been in the arsenal of biologists for almost a
century, at least since McCrady used it for quantitative determination
of Bacillus coli in water in 1915, \cite{McCrady:1915}. Since then
many studies used it to estimate the number of objects of interest in
a medium without their direct count, \cite{Fisher:1922},
\cite{Cochran:1950}, \cite{Ridout:1995}, to name a few. 

So far, several studies have produced estimates for the number of
pre-HSCs in AGM by using the dilution method which varies between just
a few to, perhaps, as many as 200, \seg \cite{Kumaravelu:2002},
\cite{Gekas:2005}, \cite{Ottersbach:2005} \cite{Bonne:2010} and
\cite{Medvinsky:2011}. Most of the reports tend to focus on the
modelling of experimental data and on the estimation methods.
However, there has been little discussions on how to design the
experiment to capture the most informative sample. Indeed, the
experiment would be spoilt, if all the mice repopulate or if all do
not. The aim of this paper is to find an optimal design of the
dilution experiment to estimate the mean number of pre-HSCs. What
distinguishes the experiments we are dealing with in this paper from
other dilution experiments extensively covered in statistical
literature are the following two features. First, it is the \emph{volume
contraint} imposed by the limited size of available substrate from
AGMs. Second, a significant time delay between the dose injection and
the result of its action prevents from using \emph{multi-stage designs} when
the further stages of experiments are based on the outcomes of the previous
one(s). To address these specific issues, we employ recently developed
methods of \textit{constrained optimisation} of functionals of
measures and the corresponding steepest descent type algorithms for
computations.  It should be stressed that our methodology is generic
in the sense that it can be applied to other dilution experiments and
not only in the stem cells research. Moreover, additional contraints
can be incorporated into the model which would allow, for instance, to
take into account the cost of mice or other materials used in the
experiment.

This study is organised in the following
way. Section~\ref{sec:Dilution-method} introduces the dilution
experiment we are dealing with, description of the corresponding
statistical model followed by its assumptions, and at the end the
optimality criterion
functions. Section~\ref{sec:Optimisation-of-functionals-of-measure}
lays out the theoretical basis for the optimisation methods we are
using given that the goal functions are represented as functions of a
\textit{measure}. Consecutive sections provide the optimal design of
the dilution experiments under various conditions and various goal
functions: in non-Bayesian setting in
Section~\ref{subsec:Optimisation-of-the-Fisher-information-w.r.t-the-measure},
then under Uniform and Gamma Bayesian priors in
Section~\ref{subsec:bayes}, and finally in Section~\ref{subsec:cost}
we represent an optimal design that integrates cost to the criterion
functions. We conclude by discussion of our findings and their
extensions in Section~\ref{sec:discussion}.

\section{Dilution experiment, statistical model and optimality criteria}
\label{sec:Dilution-method}

\paragraph{Description of experiment}
\label{sec:descr-exper}

The dilution experiment on estimation of the number of HSCs we address
here involved the AGM region of an 11 days old mouse embryo. More
exactly, in order to make a study more representable and not depending
on features of a particular embryo, an engineered AGM is made from
several such embryos. A substrate of volume $V$ is then made from this
engineered AGM and the content is thoroughly mixed.  Next, $n$ doses
containing proportions $x_1\thru x_n$ of the whole $V$ are extracted
from all or part of the substrate and left for a few days so that
pre-HSCs in these doses, if any, mature into HSCs. Finally, these $n$
doses\footnote{These are further diluted to a standard volume, but
  this, obviously, does not change the number of HSCs present
  before the dilution.} are injected into $n$ irradiated mice (the
number $n$ of mice was 30 in this experiment) so that each mouse
receives its own dose and the mice are put to rest for a few weeks for
the doses to take effect. After this, the number of repopulated mice,
which are the ones having received a doze with at least one pre-HSC,
is counted and the inference is drawn on the total number of pre-HSCs
initially present in the AGM. The main question to address when
designing such an experiment is doses of which volumes should be used
for available number of mice in order to get best possible quality of
the statistical estimator? Different criteria could be considered to
quantify the quality of the estimator. We will consider below most common ones based on the Fisher information.

\paragraph{Statistical model}
\label{sec:stat-model}

Given the description of the experiment above, the following
assumptions can be made. 
\begin{itemize}
\item \emph{Spatial homogeneity}: pre-HSCs were distributed uniformly
  throughout the substrate when the doses were taken. Thus there is no
  tendency for pairs or groups of pre-HSCs either to cluster or to
  reject one another. This is implied by the fact that the substrate
  was thoroughly mixed just before the doses are taken.
\item \emph{Orderliness}: the probability that there are more than one
  pre-HSCs in a small volume $dv$ of substrate has order $o(dv)$. This
  is a natural assumption given that an 11 days old AGM contains about
  300 thousand cells and only no more than 200 of these are pre-HSCs.
\item \emph{Independence}: each cell in substrate has the same (small)
  probability to turn into a pre-HSC independently of the other
  cells.
\item Only a pre-HSC can develop into a mature HSC. So that a dose
  contains at least one HSC at the time of injection to a mouse if and
  only if there was a pre-HSC present in the dose at the time of its
  extraction from the whole substrate.
\item Finally, each dose when injected into an irradiated mouse is
  certain to exhibit a positive result (repopulated mouse), whenever
  the dose contains at least one HSC.
\end{itemize}

The first three assumptions above suggest that the locations of
pre-HSCs in the substrate $V$ are given by a homogeneous Poisson point
process. This follows from the Poisson limit theorem for thinned point
processes, \seg~\cite[Sec.~11.3]{DalVJon:08}. Indeed, in every subset
of the substrate of a positive volume $x$, the number of cells turned
into pre-HSCs is well approximated by the Poisson distribution with
the parameter proportional to the mean number of pre-HSCs in the
substrate which is $\lambda x$ for some parameter $\lambda>0$. Because
of the independence assumption, these numbers are independent for disjoint
subsets. Changing the units if necessary, we assume from now on that
the volume of the whole substrate is 1. The parameter $\lambda$ is then the
unknown density of the Poisson point process which is also the mean
number of pre-HSCs in the substrate. Thus we operate with a measurable
space carrying point configurations inside a set $V\subset \R^3$ of
volume 1 (the space of finite counting measures $\omega$ on $V$ with
the minimal $\sigma$-algebra making the mappings $\omega\mapsto
\omega(B)$ measurable for all Borel $B\subseteq V$) supplied with
probability distribution $\P_\lambda$ so that $\omega$ under
$\P_\lambda$ is a homogeneous Poisson point process with density
$\lambda$.

The doses taken can now be associated with disjoint subsets $V_1\thru
V_n$ of $V$ with volumes $x_i>0$, $i=1\thru n$. The corresponding
numbers of pre-HSCs $\omega(V_i),\ i=1\thru n$ in the doses are then independent
Poisson distributed random variables with parameters $\lambda x_i$
while the total number of pre-HSCs $\omega(V)$ is Poisson distributed with
parameter $\lambda$.

In the simplest case all the doses have the same volume $0<x\leq
1/n$. The probability that a doze is sterile, i.e.\ it does not
contain a pre-HSC, is then
\begin{equation}
\label{eq:prob-ith-dose-is sterile-equal-dose}
p=\P\{\omega(V_i)=0\}=e^{-\lambda x}. 
\end{equation}
Thus the total number of non-repopulated mice follows Binomial
distribution with parameters $n, p$ and the maximum likelihood estimate for the average
number of HSC $\lambda$ is given by
\begin{equation}
\label{estimation-of-lambda-equal-dose}
\hat{\lambda}=-\frac{\log \hat{p}}{x},
\end{equation}
where $\hat{p}$ is the proportion of non-repopulated mice provided it
is not 0. However, the doses need not be necessarily all equal for an
optimal design. 

Let $\chi_i$ ($i=1,...,n$) be an indicator that a mouse, which
received the $i$th dose of volume $x_i$, has not repopulated. Thus,
$\chi_i$ is a \textit{Bernoulli} random variable
\begin{equation}
\label{eq:Bernoulli-random-variable}
{\chi_i}|{\lambda} \sim \operatorname{Bern} ({e^{-\lambda x_i}})
\end{equation}
with the parameter equal to the probability of the $i$th dose to be sterile.

Hence, the \textit{log-likelihood function} for the sequence $\bm\chi=(\chi_1,...,\chi_n)$ of
non-repopulated and repopulated mice is given by
\begin{equation}
\label{eq:log-likelihood-function-of-n-observations}
\ell(\bm\chi\cond\lambda,\mathbf{x})=-\sum_{i=1}^n \chi_i \lambda x_i
+\sum_{i=1}^n (1-\chi_i)\log(1-e^{-\lambda x_i}), 
\end{equation}
where $\mathbf{x}=(x_1,...,x_n)$. Maximisation of this expression over
$\lambda$ for an observed sample $\bm\chi$ provides a maximum
likelihood estimator (MLE) of $\lambda$ for a given design
$\mathbf{x}$. Our goal here is to determine the optimal design in
terms of the doses $\{x_i\}$, according to a suitably chosen
\textit{optimality criterion}, which we describe next.

\paragraph{Optimality criterion functions}
\label{sec:Optimality-criterion-functions}
Recall that for a statistical model which depends on a one-dimensional
parameter $\lambda$, the \textit{Fisher information} is defined as
\begin{equation}
  \label{eq:expected-Fisher-information-formula}
  I(\mathbf{x};\lambda)=-\E_\lambda 
  \biggl[\frac{\partial^2 \ell(\bm\chi|\lambda,\mathbf{x})}{\partial \lambda^2}\biggr].
\end{equation}
The expectation $\E_\lambda$ is taken here w.r.t the random vector
$\bm\chi$.  

Derived from (\ref{eq:log-likelihood-function-of-n-observations}), we
have in our case
\begin{equation}
  \label{eq:Fisher-information-function}
  I(\mathbf{x};\lambda) = \sum_{i=1}^n \frac{e^{-\lambda x_i}}{1-e^{-\lambda x_i}} x_i^2.
\end{equation}
The Fisher information measures the amount of information that an
observable sample carries about the unknown parameter, which the
likelihood function depends upon. On the other hand, it is the inverse
of MLE's variance, \seg~\cite{Everitt:2010}. Thus,
maximising the information corresponds to minimising the variance of
the MLE. Therefore, maximising the Fisher information
(\ref{eq:Fisher-information-function}) over $\mathbf{x}$ under
constraint $\sum_{i=1}^n x_i\leq 1$ is a useful design criterion.

It is typical in statistical experiment planning to describe design in
term of a probabilistic \emph{design measure}. Typically,
the design measure is atomic, so it has a form $\sum_j q_j
\delta_{x_j}$, where $\delta_x$ is the unit mass measure concentrated
on a point $\{x\}$. The design measure reflects the (asymptotic when $n\ti$)
frequency $q_j$ of occurrence of the value $x_j$ in the design,
\seg~\cite{atk:don92}. By this reason, we will also describe the
doses by a measure $\mu(dx)$ living on $(0,1]$, albeit not renormalised
to have mass 1. Namely, $\mu=\sum_j m_j\delta_{x_j}$ means the design
when a dose of volume $x_j$ is repeated $m_j$ times. Since it is
senseless for an experiment involving estimation of $\lambda$ to give
a mouse zero-doze of the substrate, we exclude the point 0 from the
design space. Obviously, we
have that the total mass constraint $\sum_{j=1}^k m_j=\int \mu(dx)=n$ and that we
cannot extract more doses than the total volume of the substrate:
$\sum_{j=1}^k m_j x_j=\int x\mu(dx)\leq 1$. All integrals here and
below are taken over $(0,1]$, unless specified differently.

Now the basic optimisation problem for the design of our experiment is
\begin{align}
  \label{eq:G1}
  G_1(\mu;\lambda) & = I(\mu;\lambda)=\int \frac{e^{-\lambda x}}{1-e^{-\lambda x}}
  x^2\,\mu(dx)\to\sup\\
  \intertext{over measures $\mu$ with support $\supp \mu\subseteq [0,1]$
    satisfying}
  & \mu(\{0\})=0; \notag\\
  & \mu((0,1])=n; \label{eq:tot-mass}\\
  & \int x\,\mu(dx)\leq 1.\label{eq:tot-volume}
\end{align}
A design measure describes \emph{asymptotic} frequencies, so $q_j n=m_j$
for a given finite $n$ may not be all integers. In this case it is
reasonable to consider the nearest measure with all $m_j\in \Z_+$ as
an approximation to the optimal solution. Or, if necessary, a choice
of the optimal measure among such measures can be done by evaluation
of the goal function at just a few closest approximations of this kind
to the optimal design measure.

\paragraph{Bayesian setting}
Sometimes, there is an additional information available on the
plausible values of the parameter $\lambda$ which is given in a form
of a \textit{prior distribution} $Q(d\lambda)$,
\seg~\cite[Sec. 2.2.2]{Mar:Rob07}. In this case the optimality
criterion involves taking the expectation of the goal functions w.r.t
the distribution $Q$.

For a single parameter, the following criterion functions are typically used to find an optimal design, see, e.g., \cite{atk:don92} or \cite{Ridout:1995}:
\begin{align}
    & G_2(\mu)=\E_Q I(\mu;\lambda) \label{eq:G_2}\\
    & G_3(\mu)=\E_Q \log I(\mu;\lambda) \label{eq:G_3}\\
    & G_4(\mu)=-\E_Q \bigl(I(\mu;\lambda)\bigr)^{-1} \label{eq:G_4}
\end{align} 
Criterion function $G_2$ and $G_3$ maximize, under the same
constraints \eqref{eq:tot-mass} and \eqref{eq:tot-volume}, the expectation of the
Fisher information and of its logarithm, respectively, as used in
e.g., in \cite{za:77} and \cite{ch:la89}. The criterion
function $G_4$ minimises the expected asymptotic variance
of the maximum likelihood estimator.  

Next section will describe a general framework of optimisation
of functionals of measures and the recently developed techniques for
solving such optimisation problems. Apart from optimal design of
statistical experiments \cite{Pukelsheim:83}, \cite{atk:don92}, these
are frequent in different subjects, like spline approximation of curves and
geometrical bodies where the measure describes the positions of spline
points 
\cite{Schneider:88}, maximisation the
area covered by random geometric objects with the distribution
determined by a measure \cite{Hall:88}, stochastic search,
where a measure determines the search strategy \cite{WZ:94}, \cite{Z:91}. 

\section{Optimisation of functionals of measures}
\label{sec:Optimisation-of-functionals-of-measure}
In this section we summarise necessary information about measures and
variational analysis on them. Further details on measure theory can be
found, \eg in~\cite{dun:sch09} or  \cite{HP:57}. 

Let $X$ be a locally compact separable topological space with the Borel
\textit{$\sigma$-algebra} $\bbb$ of its subsets. Let $\M$ ($\M_+$) denote the
set of signed (respectively, non-negative) finite measures on $\bbb$,
\ie\ countably additive functions from $\bbb$ to $\R$ ($\R_+$,
respectively). $\M$ becomes a linear space if the sum of measures and
multiplication by a number is defined by
$(\eta+\nu)(B)=\eta(B)+\nu(B)$ and $(t\eta)(B)=t\eta(B)$ for any
$t\in\R$ and $\eta,\nu\in\M$. $\M_+$ is a \emph{cone} in $\M$ since
$\mu+\nu\in \M_+$ and $t\mu\in\M_+$ whenever $\mu,\nu\in\M_+$ and
$t\geq 0$.  The \textit{support} $\supp \mu$ of a positive measure
$\mu$, is defined as the complement to the union of all open sets of
zero $\mu$-measure. Measures are \emph{orthogonal} if their supports
are disjoint. Any signed measure $\eta$ can be represented as the
difference $\eta^+-\eta^-$ of two orthogonal non-negative measures
$\eta^+,\eta^-\in\M_+$ (the Jordan decomposition). The set $\M$
becomes a \emph{Banach space} if supplied with the \textit{total
  variation norm} $\|\eta\|=\eta^+(X)+\eta^-(X)$.

Optimisation of functions defined on a Banach space, which are commonly
called \emph{functionals}, relies on the
notions of differentiability. In our case a functional $f:\ \M\mapsto
\R$ is called \emph{strongly} or \emph{Fr\'echet differentiable} at
$\nu\in\M$ if
\begin{equation}
\label{eq:Frechet-differentiable}
f(\nu+\eta)-f(\nu)=Df(\nu)[\eta]+o(\|\eta\|) \ \text{as} \ \|\eta\| \downarrow 0,
\end{equation}
where $Df(\nu)$ is a bounded linear continuous functional on $\M$
called the \emph{differential}.

When a function is strongly differentiable at $\nu$ then
there also exists a \emph{weak} (or \emph{directional} or
\textit{Gateaux}) derivative, \ie
\begin{equation}
  \label{eq:Gateaux differentiable}
  \lim_{t \downarrow 0} t^{-1} [f(\nu+t\eta)-f(\nu)]=Df(\nu)[\eta]
\end{equation}
for any `direction' $\eta\in\M$.

The differential $Df(\nu)$ often has an integral form:
\begin{displaymath}
Df(\nu)[\eta]=\int g(x;\nu)\,\eta(dx)
\end{displaymath}
for some function $g(\,\mydot\,;\nu):\ X\mapsto \R$ which is then
called the \emph{gradient function} to $f$ at $\nu$. Not all linear
functionals are integrals, unless the space $X$ is a finite set in
which case $\M$ can just be identified with an Euclidean space. In
most applications, however, including the experimental design,
differentiable functionals do possess a gradient function, so this
assumption is not too restrictive in practice.


In this study we are interested in optimisation of functionals of
positive measures. A general constrained optimisation problem can be
written as follows:
\begin{equation}
\label{eq:optimisation-on-M-with-out-constrain-1}
f(\mu)\longrightarrow \inf, \quad  \mu \in \A,
\end{equation}
where $\A\subseteq\M_+$ is a set of measures describing the
constraints. If $f$ is strongly differentiable then the first order
necessary optimality condition states that if $\mu^*$ provides a local
minimum in the
problem~\eqref{eq:optimisation-on-M-with-out-constrain-1} then
\begin{equation}
\label{eq:optimisation-on-M-with-out-constrain-2}
Df(\mu^*)[\eta]\geq 0 \quad\text{for all} \ \eta\in T_\A(\mu^*),
\end{equation}
where 
\begin{equation}
\label{eq:tangent-cone-to-A}
T_\A(\mu)=\liminf_{t \downarrow 0}\frac{\A-\mu}{t}
\end{equation}
is the \textit{tangent cone} to $\A$ at $\mu$. Here the `$+$'
(respectively `$-$') operation on sets indicates all pairwise sums of
(respectively difference between) the points from the corresponding
sets. The tangent cone is the closure of all admissible directions
$\eta\in\M$ at $\mu$ meaning that $\mu+t\eta\in\A$ for all
sufficiently small $t>0$, \seg~\cite{C:90}. In other words, derivative
in all admissible directions should be non-negative at a point of
local minimum.  For any $\A$ of interest, one needs to characterise
the tangent cone $T_\A(\mu)$.

General optimisation theory for functionals of measures has been
developed in a series of papers \cite{MZ:2000_2}, \cite{MZ:2000_1} and
\cite{mz:04}. For us here an optimisation with finite number of equality
and inequality constraints is relevant.

Consider the following optimisation problem:
\begin{align}
\label{eq:Optimisation-with-equality-and-inequality-constraints-on-M_+}
 & f(\mu)\rightarrow \inf ,\ \mu\in \M_+\ \text{subject to}\\
& \begin{cases} 
  H_i(m)=0\quad i=1,...,k,\ k\leq d\\
  H_j(m)\leq 0 \quad  j=k+1,...,d. \label{eq:constraints}
\end{cases}
\end{align}
where $f:\ \M_+\mapsto \R$ and $H=(H_1\thru H_d):\ \M_+\mapsto \R^d$
are strongly differentiable functions. Alternatively, the
constraints~\eqref{eq:constraints} can be written in the form
$H(\mu)\in C$, where $C\subset \R^d$ is the cone $\{\mathbf{y}\in\R^d:
y_1=\dots=y_k=0,\ y_{k+1}\leq 0\thru y_d\leq 0\}$.

\begin{Dfn}[\cite{Ro:76}]
\label{Def:Robinson's-regularity-condition}
A measure $\mu$ is called \emph{regular} for optimisation problem
(\ref{eq:Optimisation-with-equality-and-inequality-constraints-on-M_+})
if the origin $0$ of $\R^d$ belongs to the interior of the set
\begin{equation}
\label{eq:Robinson's-regularity-condition}
H(\mu)-C+DH(\mu)[\M_+-\mu]\subseteq\R^d.
\end{equation}
\end{Dfn} 
Robinson's regularity condition which as shown in \cite{ZK:79},
guarantees the existence and boundedness of the \textit{Kuhn–Tucker
  vector} appearing in the next theorem. See also \cite{MZ:79} for
the discussion of different forms of regularity condition.

The following theorem gives the first-order necessary conditions for a
minimum in the
problem~(\ref{eq:Optimisation-with-equality-and-inequality-constraints-on-M_+}).

\begin{Th} \cite[Th. 4.1]{MZ:2000_2}.
\label{Th:General-equivalence-theorem}
Let $\mu^*$ be a regular local minimum of $f$ over $\mathbb{M}_+$,
subject to \eqref{eq:constraints}. Assume that $f$ and $H$ are
Fr\'echet differentiable at $\mu^*$, and there exist the corresponding
gradient functions $g(x,\mu^*)$ and $h_i(x,\mu^*), i=1,...,d$. Then
there exists Kuhn–Tucker vector $(u_1\thru u_d)\in\R^d$ such that
$u_j\leq 0$ (resp. $u_j=0$) for those $j\in \{k+1,...,d\}$ satisfying
$H_j(\mu^*)=0$ (resp. $H_j(\mu^*)< 0$), such that
\begin{equation} \label{eq:Lagrange-multiplier-General-equivalence-theorem}
  g(x,\mu^*) \begin{cases}
    =\sum_{i=1}^d u_i h_i(x,\mu^*) & \quad\text{for $\mu^*$-almost all}\ x,\\
    \geq \sum_{i=1}^d u_i h_i(x,\mu^*) &\quad\text{for all}\ x\in X.
  \end{cases}
\end{equation}
\end{Th}

One can show that in the case of finitely many constraints
\eqref{eq:constraints} satisfying
\eqref{eq:Lagrange-multiplier-General-equivalence-theorem}, the regularity condition
\eqref{eq:Robinson's-regularity-condition} becomes the so-called
\textit{Mangasarian--Fromowitz constraints qualification} (see
\cite[p. 274]{C:90}), that is the linear independence of the gradient functions
$h_1(\mydot,\mu^*)\thru h_k(\mydot,\mu^*)$ and the existence of a measure $\eta\in \M$ such that
\begin{align} 
\int h_i(x,\mu^*)\,\eta(dx)=0 & \quad \text{for all}\ i=1,...,k;
 \label{eq:Mangasarian-Fromowitz-constraints-qualification-1} \\
\int h_j(x,\mu^*)\,\eta(dx)< 0 & \quad \text{for all}\ j\in \{k+1,...,d\}\notag\\ 
& \quad \text{verifying $H_j(\mu^*)=0$}.
\label{eq:Mangasarian-Fromowitz-constraints-qualification-2}
\end{align}
Without the inequality constraints, condition
\eqref{eq:Mangasarian-Fromowitz-constraints-qualification-2},
trivially holds for $\eta$ being the zero measure.

The design problems we consider naturally fall in the above described
general framework of optimisation of functionals defined on finite
measures. Theorem~\ref{Th:General-equivalence-theorem} provides the
necessary conditions for optimality of a design. Moreover, it allows
one to easily incorporate into the model other constraints on the
optimal design measure, if needed. Constraints~\eqref{eq:tot-mass} and
\eqref{eq:tot-volume} correspond to linear functionals $H_1(\mu)=\int
\mu(dx)-n$ and $H_2(\mu)=\int x\,\mu(dx)-1$ with the corresponding
gradient functions $h_1(x;\mu)\equiv 1$ and $h_2(x;\mu)=x$. These
constraints are regular for any $\mu$ since
\eqref{eq:Mangasarian-Fromowitz-constraints-qualification-1} and
\eqref{eq:Mangasarian-Fromowitz-constraints-qualification-2} are
satisfied, for instance, for a measure $\eta=\delta_0-\delta_1$. We
therefore have the following important corollary of
Theorem~\ref{Th:General-equivalence-theorem} which we use in the next
section. Note that we mostly \emph{maximise} the goal function so that
the inequalities in
\eqref{eq:Lagrange-multiplier-General-equivalence-theorem} change to
opposite.

\begin{Th}
\label{th:optimisation-with-a-limited-cost}
Let $\mu^*$ be a local maximum of a strongly differentiable function $f:\ \M_+\mapsto \R$ possessing a gradient function
$g(x;\mu^*)$, subject to contraints \eqref{eq:tot-mass} and
\eqref{eq:tot-volume}. Then, there exist constants $u_1$ and $u_2$,
where $u_2\geq 0$ if $\int x\mu^*(dx)=1$ and $u_2=0$ if $\int
x\mu^*(dx)<1$, such that
\begin{equation}\label{eq:main}
g(x,\mu^*)\begin{cases} 
& =u_1+u_2 x \quad \mu^*\text{-almost everywhere},
 \\
& \leq u_1+u_2 x \quad \text{for all}\ x \in X.
\end{cases}
\end{equation}
\end{Th}

\section{Construction of optimal design}
In this section we apply the necessary condition for extremum of a
functional of measures to find optimal designs for a range of goal
functions and most common prior distributions in the Bayesian
settings. First we assume that the parameter $\lambda$, the mean
number of HSCs in the substrate, is known from previous experiments, and
obtain the optimal design, in terms of maximisation of the Fisher
information, for each $\lambda$. 

\subsection{Optimal design for a fixed average number of HSCs}
\label{subsec:Optimisation-of-the-Fisher-information-w.r.t-the-measure}

\begin{figure}[ht]
(a)~\includegraphics[width=6cm]{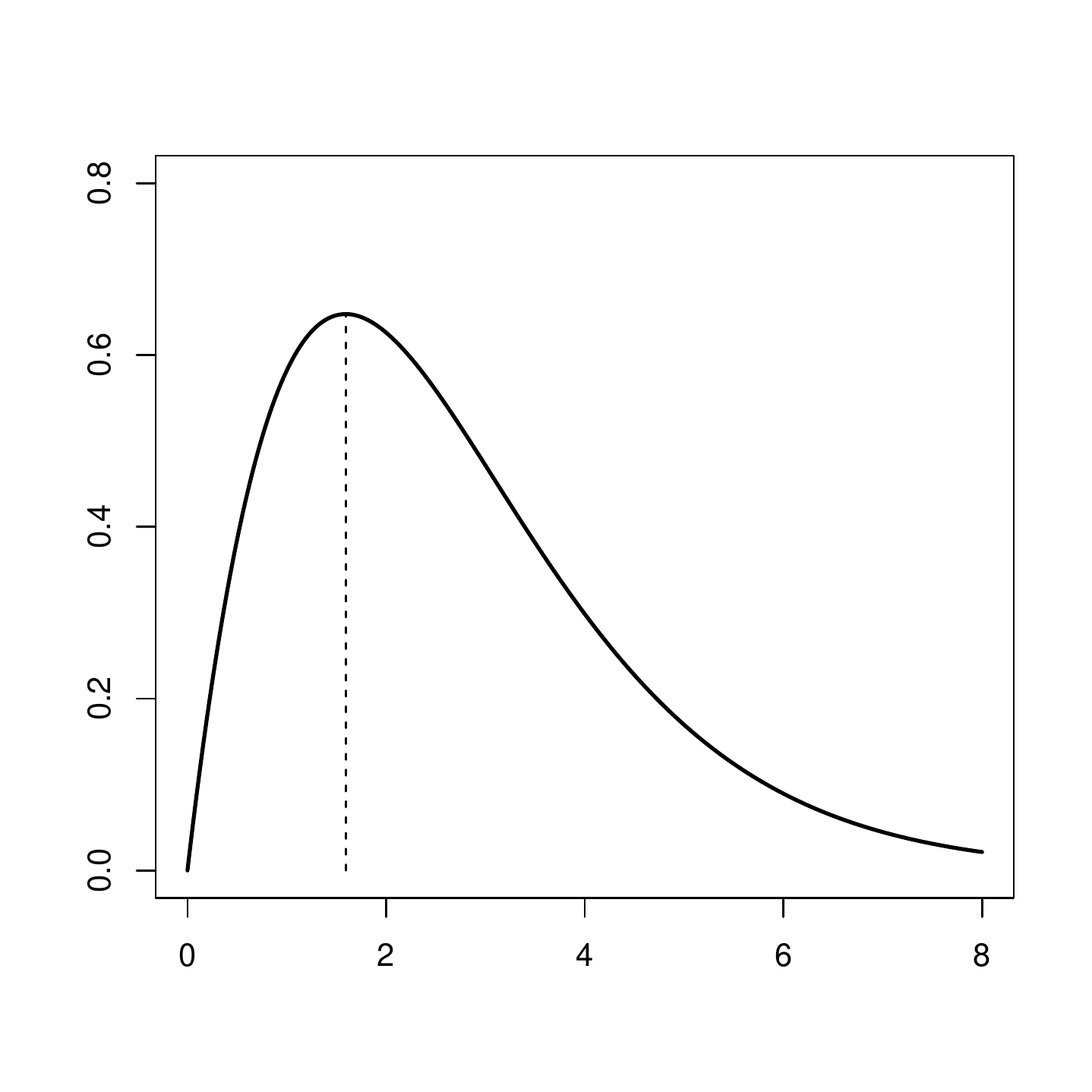}
(b)~\includegraphics[width=6cm]{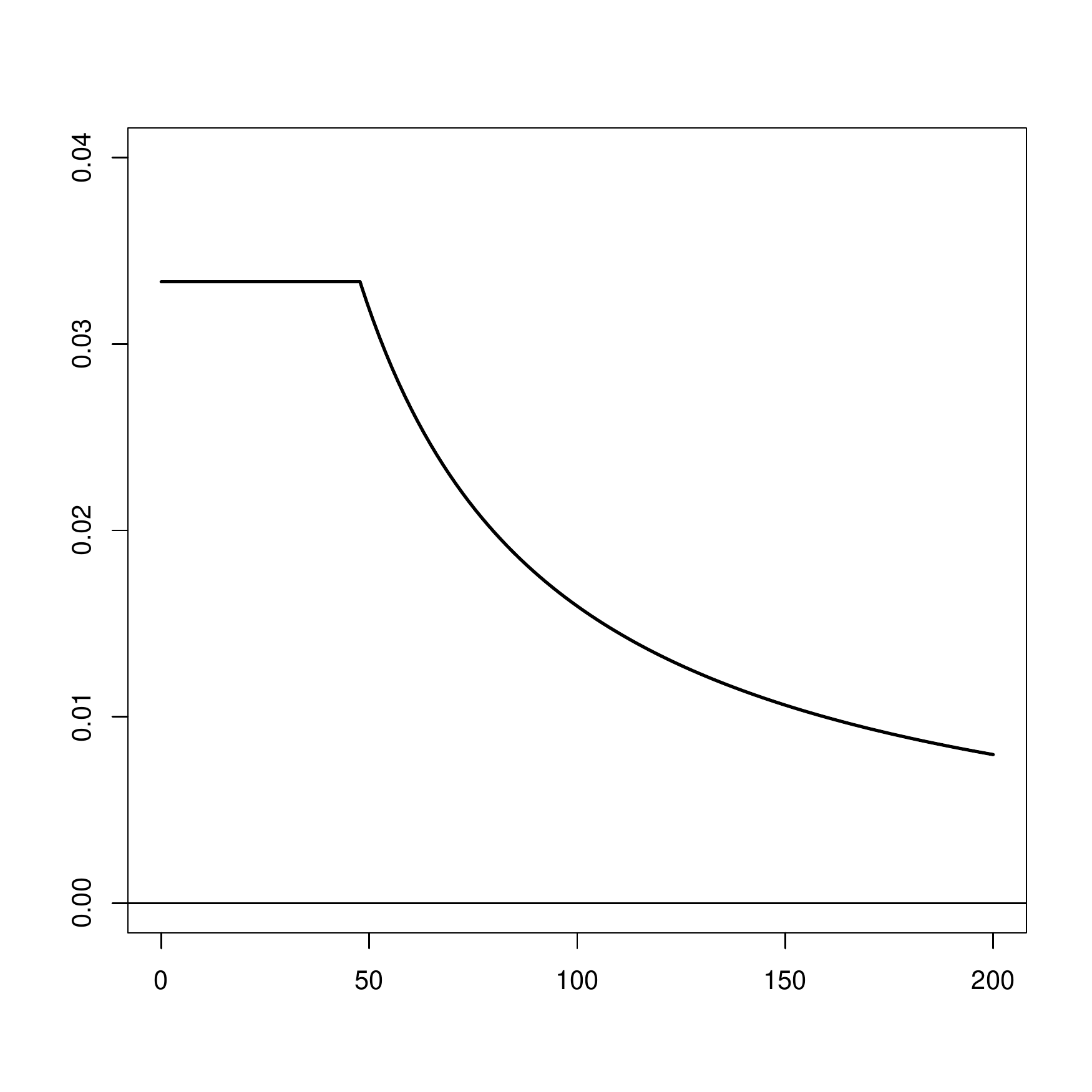}
(c)~\includegraphics[width=6cm]{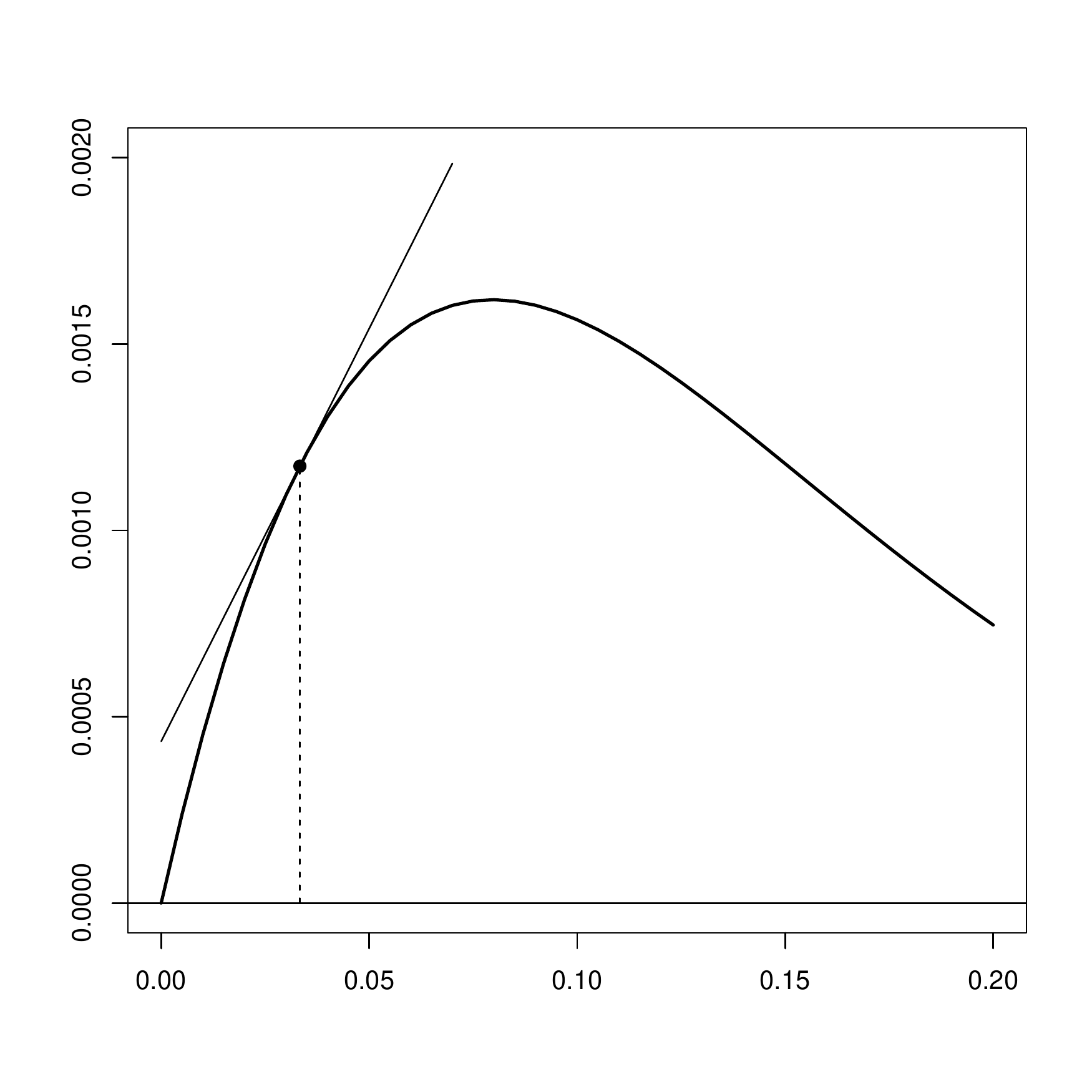}
(d)~\includegraphics[width=6cm]{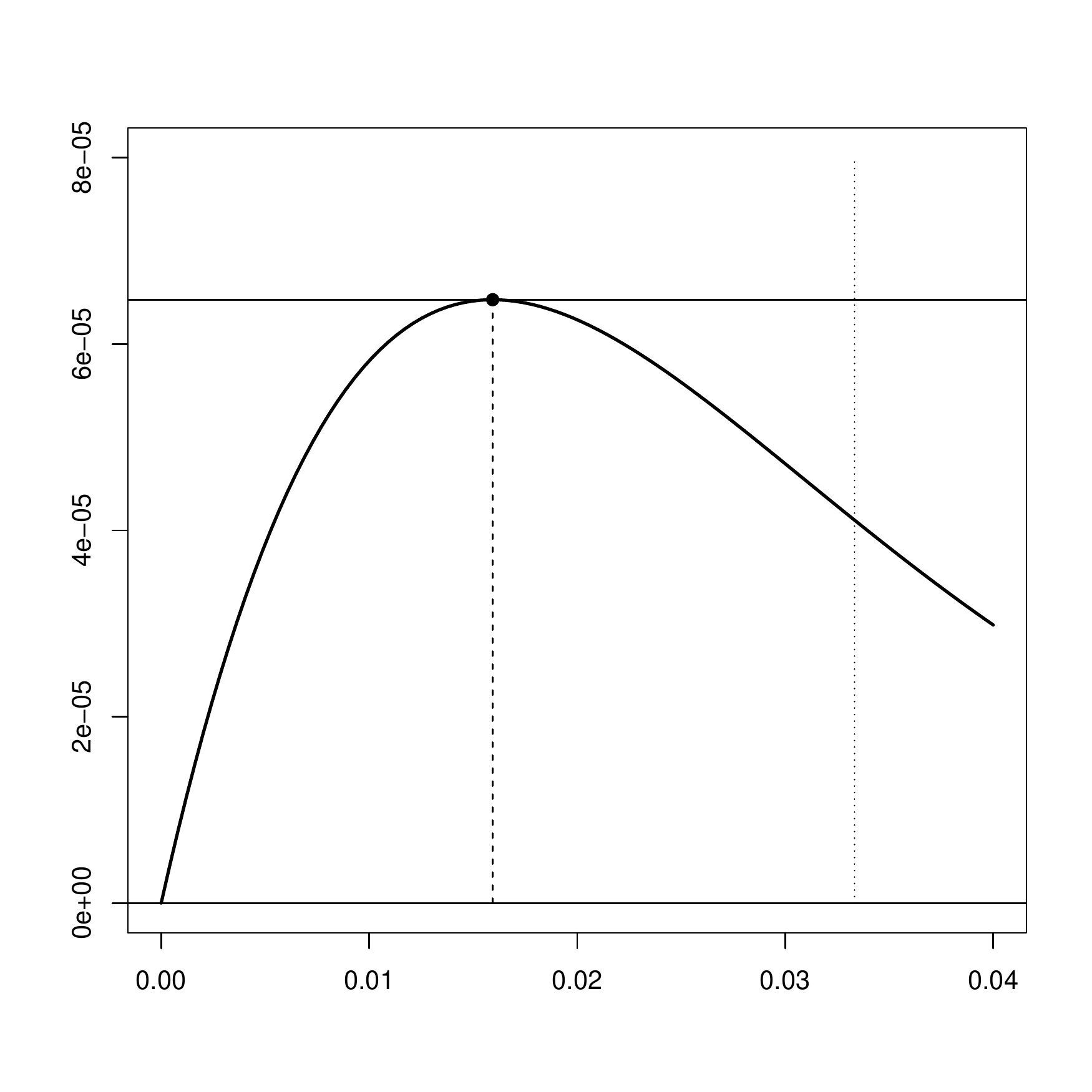}
\caption{Plot~(a): function $r$, the maximum is attained at the
  point $y_{max}\approx 1.59362$. (b): the optimal dose volume for
  $n=30$ mice as a function of $\lambda$: $1/30$ for $\lambda\leq 47.8$
  and $1.59362/\lambda$ otherwise. The line $u_1+u_2x$ and the gradient function $g_1(x;\lambda)$ satisfying
  conditions \eqref{eq:main} for $\lambda=20$ (Plot~(c), the tangent line at
  point 1/30) and for $\lambda=100$ (Plot~(d), the tangent line at the point of
  maximum $y_{max}/\lambda<1/30$).
\label{fig:r}}
\end{figure}

Here we are dealing with optimisation problem~\eqref{eq:G1} under
constraints \eqref{eq:tot-mass} and
\eqref{eq:tot-volume}. The goal function is a linear function of
$\mu$, so that its differential is the function itself with the
gradient function
\begin{equation}
\label{eq:gradient-of-G_1(m)}
g_1(x;\lambda)= \frac{e^{-\lambda x}}{1-e^{-\lambda x}} x^2,
\end{equation}
independent of $\mu$. Note that $g_1(x;\lambda)=\lambda^{-2}
r(\lambda x)$, where 
\begin{equation}\label{eq:def-r}
  r(y)= \frac{e^{-y}}{1-e^{-y}} y^2.
\end{equation}
The graph of $r$ is shown on Figure~\ref{fig:r}. It attains its unique
maximum at point $y_{max}\approx 1.59362$ and it is strictly concave on
$[0,y_{max}]$. 

The gradient function $g_1$ is just the function $r$ scaled along both axes,
and it attains its maximum at the point $y_{max}/\lambda$.
It follows from Theorem~\ref{th:optimisation-with-a-limited-cost},
that if $\mu^*$ is a measure at which $G_1$ attains its maximum, then 
$g_1(x;\lambda)\leq u_1+u_2x$ for all $x\in[0,1]$ with $u_2\geq0$. Moreover,
$g_1(x;\lambda)=u_1+u_2x$ for $x\in\supp\mu^*$. But this is only possible if
$u_1+u_2x$ is a tangent line to $g_1$ at some $x^*\in[0,y_{max}/\lambda]$ and hence
the support of optimal $\mu^*$ consists of only this point
$x^*$. Using \eqref{eq:tot-mass} and substituting $\mu=n\delta_{x}$ into \eqref{eq:G1} and
\eqref{eq:tot-volume}, we come to the optimisation problem of one variable:
\begin{displaymath}
  \label{eq:1}
  r(\lambda x)\rightarrow \sup \quad\text{over $x\in[0,1/n]$},
\end{displaymath}
so that $x^*=1/n$ for $\lambda\leq y_{max}n$ and $x^*=y_{max}/\lambda$, otherwise.
Thus we have proved the following theorem.

\begin{Th}
  The optimal design for the problem \eqref{eq:G1} under constraints
  \eqref{eq:tot-mass} and \eqref{eq:tot-volume} consists in $n$ equal
  doses of volume $1/n$ for $\lambda \leq y_{max}n$ and of volume
  $y_{max}/\lambda$ for $\lambda > y_{max}n$, where $y_{max}\approx 1.59362$ is
  the maximum point of the function $r$ \eqref{eq:def-r}.
\end{Th}

This indicates that for those $\lambda > 1.59362n$ we need to sample a
proportion of the substrate $1.59362n/\lambda$, and for those $\lambda \leq
1.59362n$ we have to take all the substrate. Therefore, if a good
prior point estimate of $\lambda$ is available, a near optimal doses
of volume $x^*$ can be selected, see Figure~\ref{fig:r}.

\subsection{Optimal design with prior distribution on $\lambda$}
\label{subsec:bayes}

Typically researchers already have an idea on what are the most likely
values of $\lambda$. This can be suggested by previous experiments or by
analogy with other similar cases and it is given in the form of a
\emph{prior distribution} $Q(d\lambda)$. The optimal experimental
design now becomes dependent not only on the choice of a criterion, but
also on the parameters of $Q$. Formally, the locations
of pre-HSCs in the substrate now conform to a Cox (or mixed Poisson)
point processes with parameter measure $Q$.

In this section we carry out optimisation of the three goal functions
introduced above in~\eqref{eq:G_2}--\eqref{eq:G_4} for the most
common prior distributions: Uniform and Gamma.

\paragraph {Uniform prior distribution.} The Uniform prior distribution
is appropriate when there is no knowledge on the mean number  $\lambda$
of pre-HSCs in substrate apart from its range. Certainly, $\lambda$ should
be greater than 1 since $pre-HSCs$ do develop in the AGM. So, it is
reasonable to assume that $\lambda\sim\Unif(1,u)$, where $u>1$ is a
known upper bound. It was already mentioned that the number of pre-HSCs
hardly exceeds 200, so one can set $u=200$, or, keeping in mind the
variance of the Poisson distribution, $u=170$, for instance.

\begin{figure}[ht]
  \centering
  \includegraphics[width=7cm]{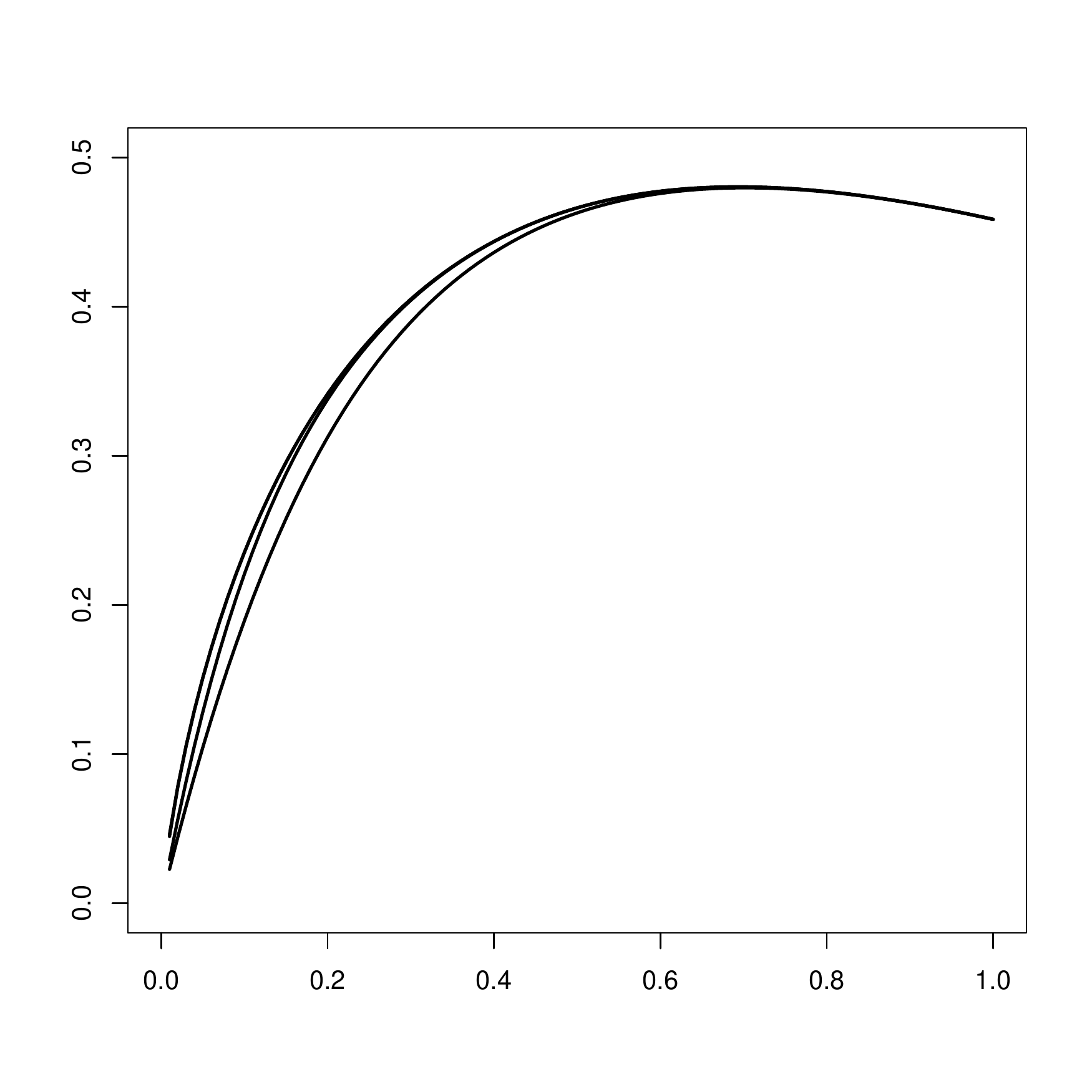}
  \caption{Function $(u-1)g_2(x;u)$ for the values $u=10,20,200,1000$ (in the order of increase). The last two curves are almost indistinguishable.}
  \label{fig:g2-unif}
\end{figure}

The goal function~\eqref{eq:G_2} now becomes
\begin{displaymath}
\label{eq:G2(m) with a unifrom prior distr}
G_2(\mu ;u)=\frac{1}{u-1}\int_{1}^{u}\int_{(0,1]}\frac{e^{-\lambda
    x}}{1-e^{-\lambda x}}x^2 \mu(dx)\,d\lambda.
\end{displaymath}
By \textit{Fubini's theorem} we can change the order of integrals
above to arrive at
\begin{equation}
\label{eq:G2(m) simplified version}
G_2(\mu ;u)=\frac{1}{u-1}\int_{(0,1]} x\log\frac{1-e^{-ux}}{1-e^{-x}}\,\mu(dx).
\end{equation}
Thus $G_2(\mu ;u)$ is a linear functional of $\mu$ with the gradient function
\begin{equation}
\label{eq:gradient function of G_2(m) w.r.t m}
g_2(x;u)=\frac{x}{u-1}\log\frac{1-e^{-ux}}{1-e^{-x}}.
\end{equation}
This function varies very little for all practically interesting values
of $u\geq 5$ with the maximum value attained around 0.69---0.7, see
Figure~\ref{fig:g2-unif}. Reasoning the same way as we did in
Section~\ref{subsec:Optimisation-of-the-Fisher-information-w.r.t-the-measure},
we can conclude that for any $u\geq 2$ the optimal design is attained
on one atom design measure concentrated on the point $1/n$ (\cf
Plot~(c) in Figure~\ref{fig:r}). So the whole volume of the substrate
should be divided into $n$ equal doses for this goal function.

Consider now optimisation criterion \eqref{eq:G_3} which now takes the form

\begin{equation}
\label{eq:G3(m) with a Uniform prior}
G_3(\mu ;u)=\frac{1}{u-1}\int_{1}^{u}\log\biggl(\int_{(0,1]} \frac{
  e^{-\lambda x}}{1-e^{-\lambda x}}x^2 \mu(dx)\biggr)d\lambda. 
\end{equation}
This goal function is not linear, but still Fr\'echet differentiable
with the gradient function given by
\begin{equation}
\label{eq:gradient function of G_3(m) w.r.t m}
g_3(x,\mu ;u)=\frac{1}{u-1} \int_{1}^{u} I ^{-1}(\mu; \lambda)\frac{e^{-\lambda x}}{1-
  e^{-\lambda x}}  x^2 \, d\lambda,
\end{equation}
where
\begin{equation}\label{eq:Idef}
  I(\mu;\lambda)=\int \frac{e^{-\lambda x}}{1-
    e^{-\lambda x}}  x^2 \,\mu(dx)
\end{equation}
is the Fisher information written in terms of a design measure.
Since the gradient now depends on $\mu$, numeric methods have to be
employed to find the optimal design for a given value of the upper
bound $u$ and the number of mice $n$. This was done by means of
R-library \texttt{medea} which finds an optimal solution to a measure
optimisation problem under linear constrains of equality
type, see \cite{MZ:2002}. Since \texttt{medea} does not yet allow for inequality
constraints, to deal with~\eqref{eq:tot-volume}, the procedure looks for an optimum
measure for a given value $b$ of the integral $\int x\mu(dx)$ and then
optimises over $b\leq 1$. The R-code is freely available from one of
the authors' web-page\footnote{http://www.math.chalmers.se/\mytilde sergei}.

\begin{figure}[ht]
  \includegraphics[width=6.5cm,height=10cm]{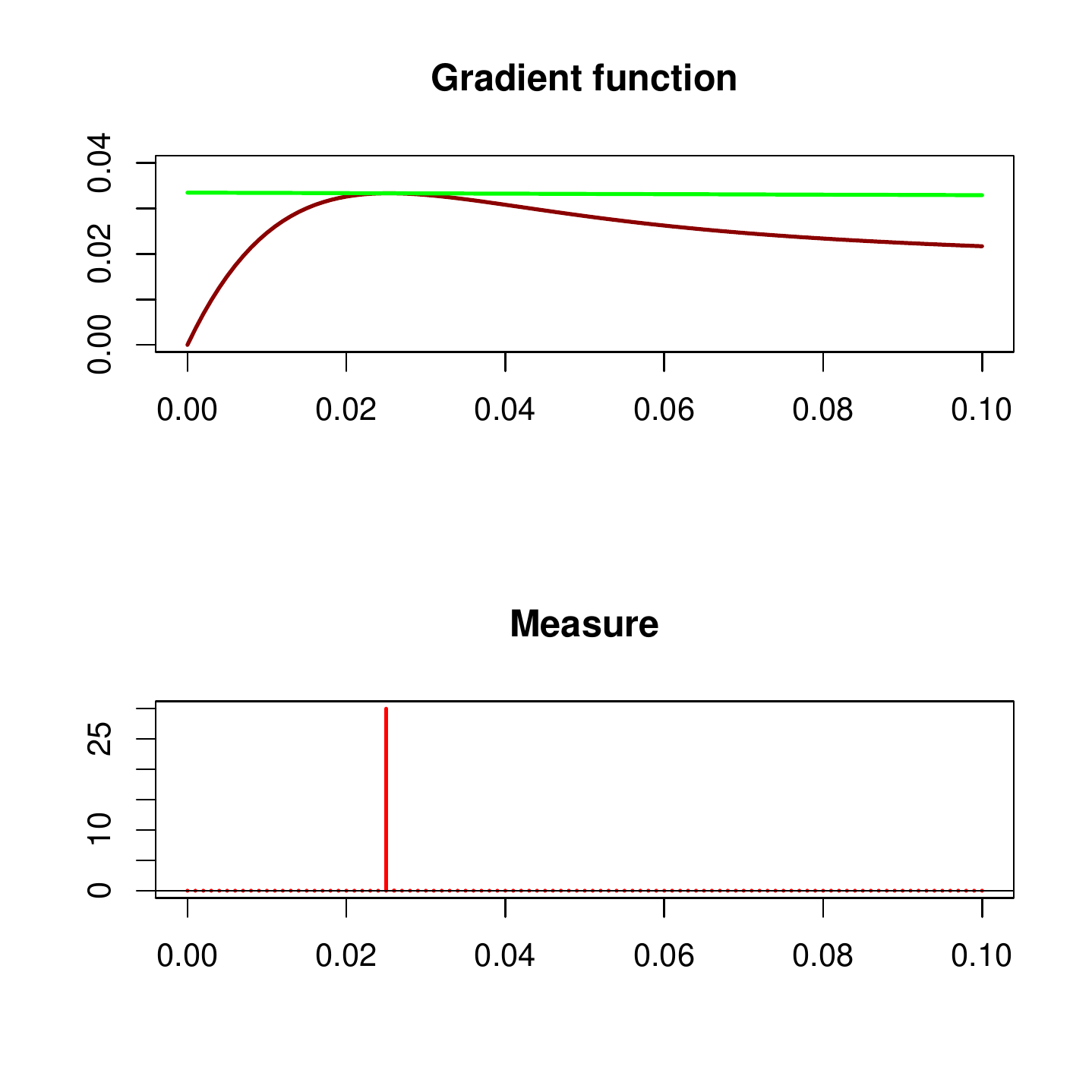}\hspace*{-5mm}
  \includegraphics[width=6.5cm,height=10cm]{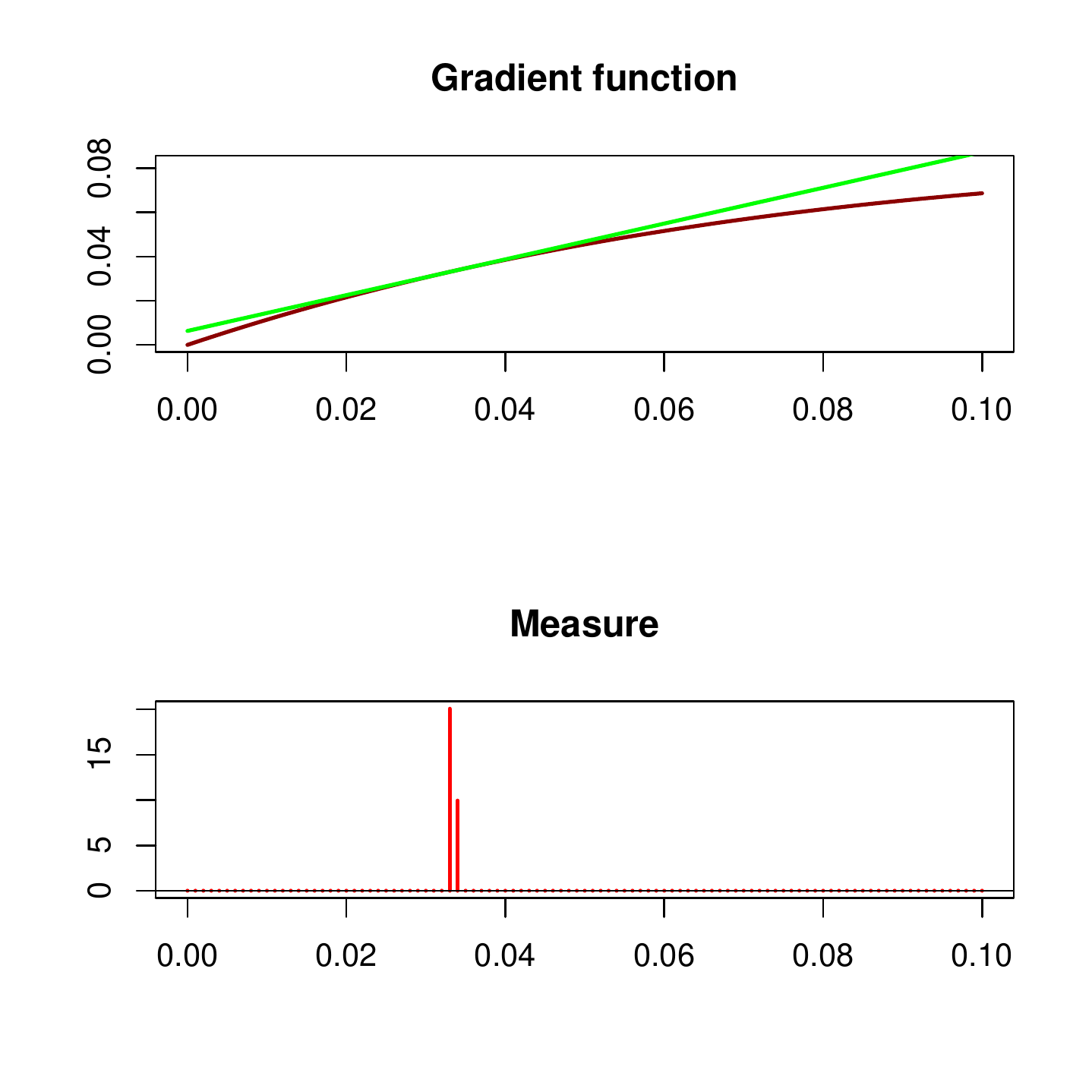}
  \caption{Numeric solution obtained by \texttt{medea} for $u=120$
    (left pair of plots) and for $u=20$ (right pair) for $n=30$ mice. 
    For $u=120$, the one-atom optimal measure is $30\,\delta_{x^*}$, 
    where $x^*=0.02522$. The numeric solution is two-atomic: 
    $29.97\,\delta_{0.025}+0.03\,\delta_{0.026}$. So, the total mass
    of $30$ was distributed between the grid points surrounding the 
    true atom position. Similarly, for $u=20$ the one-atom optimal 
    measure is $30\,\delta_{x^*}$ with $x^*=1/30$. The numeric
    solution is a two-atom measure $20.07\,\delta_{0.033}+9.93\,\delta_{0.034}$. 
    \label{fig:g3-unif-medea}}
\end{figure}

Numeric experiments conducted for various values of $u$ show that the
obtained optimal solution is always one-atom as in the previous cases.
Typically, the numeric solution gives two atoms at the neighbouring
points of the discretised space $[0,1]$ for $x$ which indicates a
single atom is situated in between these grid points, see
Figure~\ref{fig:g3-unif-medea}. Although we cannot formally prove
that the optimum design is one-atom, such designs (\ie equal doses)
are certainly of an interest.

\begin{figure}
\includegraphics[width=6cm]{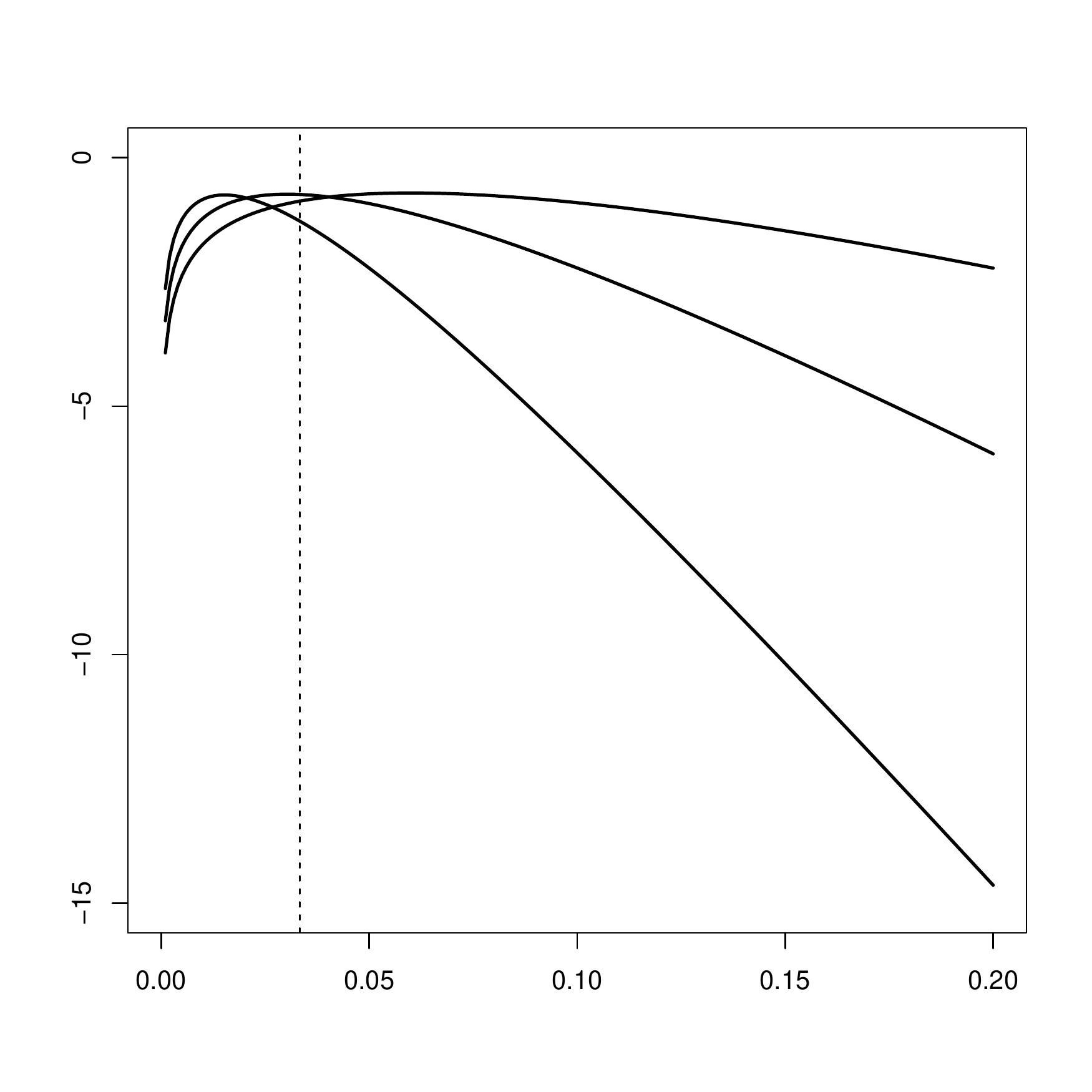}
\includegraphics[width=6cm]{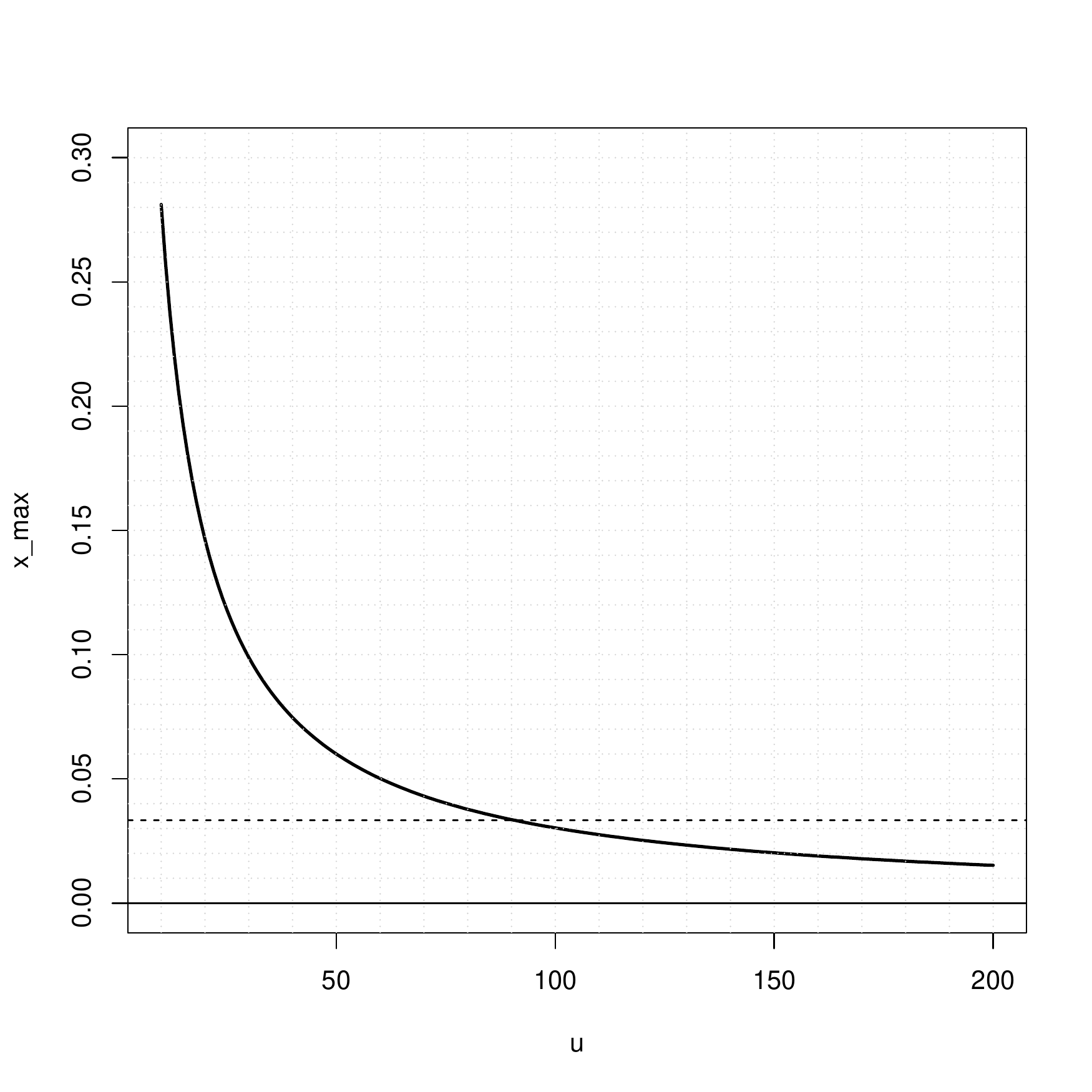}
\caption{Left plot: function $\tilde{G}_3$
  for $u=50, 100$ and 200. As $u$ grows, the point of maximum
  $x_{max}$ approaches 0. The vertical dashed line is through the
  point $1/n$ with $n=30$.  Right plot: $x_{max}$ as a function of $u$ (horizontal dashed
  line is through $1/n$ for $n=30$).\label{fig:G3tilde}}
\end{figure}

One-atom measures under constraint~\eqref{eq:tot-mass} have a form
$\mu=n\delta_{x}$ for some $x\in(0,1/n]$. Making use of
\eqref{eq:def-r}, we come to a one variable optimisation problem:
maximise
\begin{equation}
\label{eq:G3(x1)}
G_3(n\delta_{x})=\log n-2\E_Q\log\lambda+\E_Q\log r(\lambda x)
\end{equation}
subject to $x^*\leq1/n$ and $Q$ is $\Unif[1,u]$. Only the last term
depends on $x$, so the equivalent problem is to maximise 
\begin{displaymath}
  \tilde{G}_3(x;u)=\E_Q\log r(\lambda x)=\frac{1}{u-1}\int_1^u \log r(\lambda x)\,d\lambda.
\end{displaymath}
As $u$ grows, the maximum $x_{max}$ of this function, which can easily
be computed numerically, approaches 0, see the right plot on
Figure~\ref{fig:G3tilde}. So when $x_{max}\geq 1/n$, the constrained
maximum is attained at the point $1/n$. Otherwise, for large $u$,
$x_{max}<1/n$ and the solution is to take $n$ doses of volume
$x_{max}$. This is exemplified at Figure~\ref{fig:G3tilde} for $n=30$:
the optimal dose is given by
\begin{equation}\label{eq:G3n30}
x^*=
 \begin{cases}
 \frac{1}{n} & \quad u\leq u^*;\\
 x_{max}  & \quad u>u^*.  
 \end{cases}
 \end{equation}
 where $u^*\approx 90.66$.  This indicates that for those $u\leq u^*$
 one needs to take all the substrate to make $n$ equal doses, and the volume
 $nx_{max}$ otherwise.

\begin{figure}
\includegraphics[width=6cm]{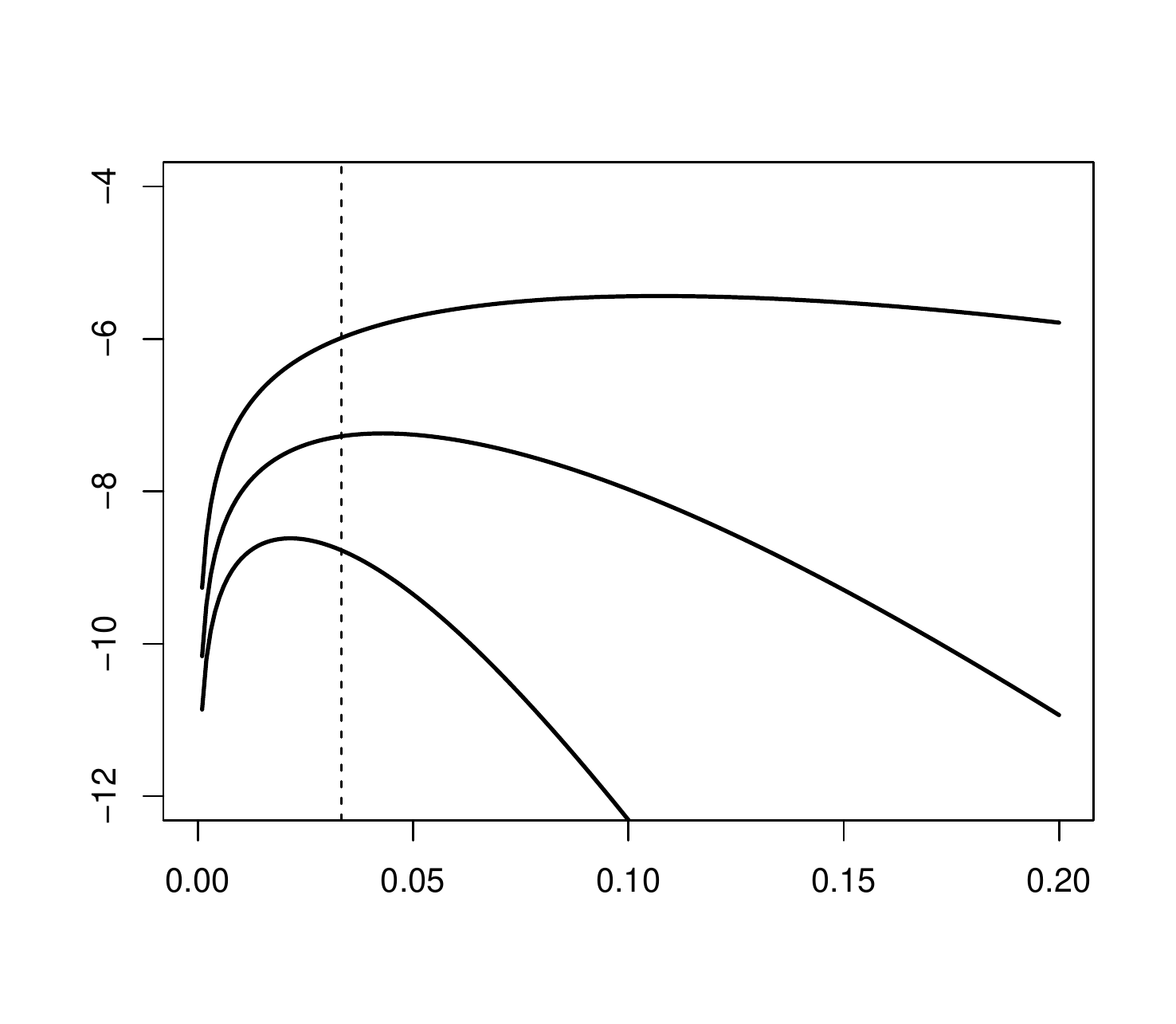}
\includegraphics[width=6cm]{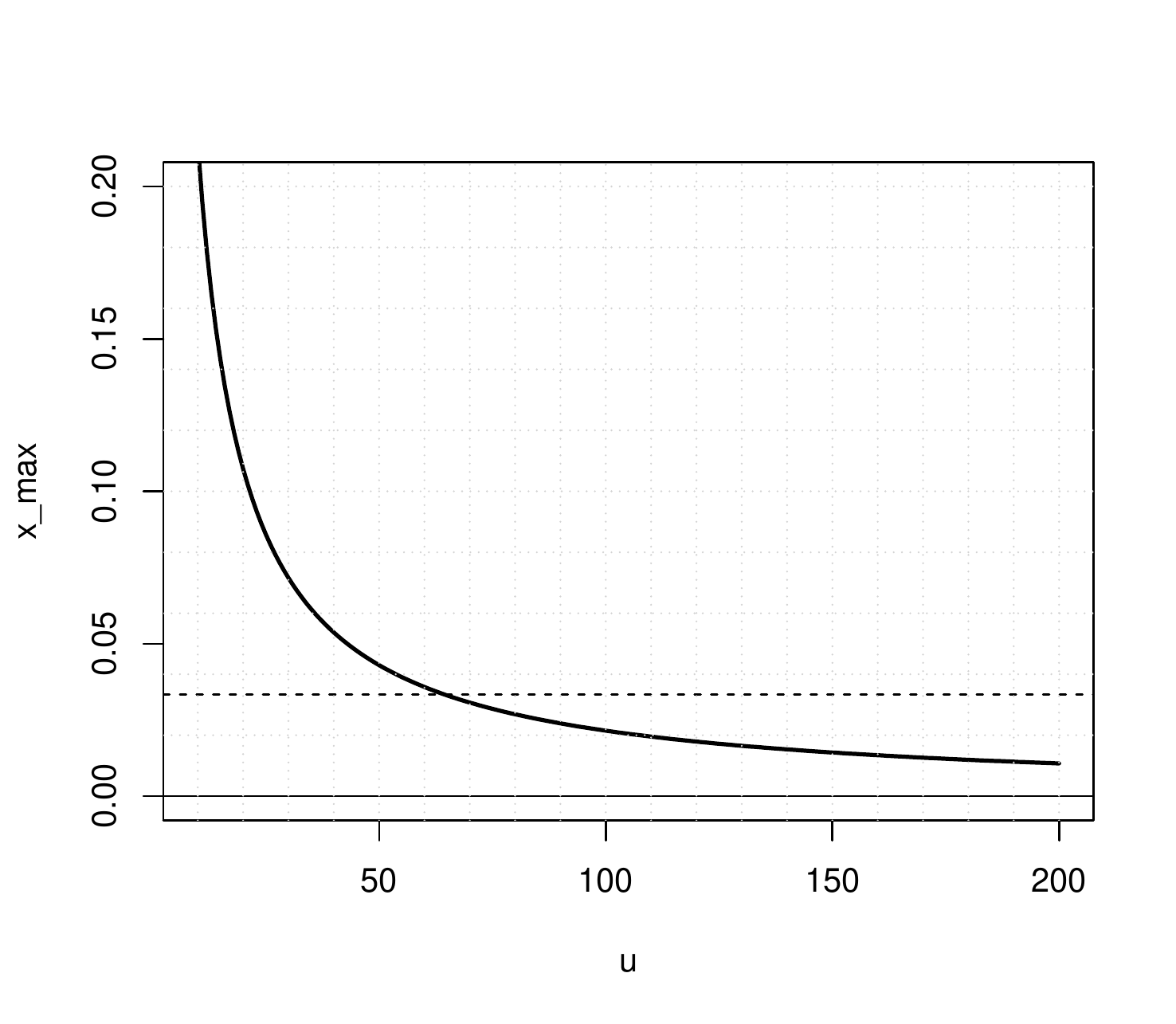}
\caption{Left plot: function $-\log(|nG_4(n\delta_{x};u)|)$
  for $u=20, 50$ and 100 (from upper to lower curve). As $u$ grows, the point of maximum
  $x_{max}$ approaches 0. The vertical dashed line is through the
  point $1/n$ with $n=30$.  Right plot: $x_{max}$ as a function of $u$ (horizontal dashed
  line is through $1/n$ for $n=30$).
\label{fig:G4}}
\end{figure}

Finally, consider objective function $G_4$ given by~\eqref{eq:G_4}
which corresponds to the expected asymptotic variance of the maximum
likelihood estimator.  In the case of a Uniform prior
$\lambda\sim\Unif(1,u)$,
\begin{equation}
\label{eq:G4(m) with a Uniform prior}
G_4(\mu ;u)=-\frac{1}{u-1}\int_1^u I^{-1}(\mu;\lambda)\,d\lambda,
\end{equation}
where $I(\mu;\lambda)$ is given by~\eqref{eq:Idef}. Again, $G_4$ is
Fr\'echet differentiable with the gradient function
\begin{equation}
\label{eq:gradient function of G_4(m) w.r.t m}
g_4(x,\mu ;u)=\frac{1}{u-1}\int_{1}^{u} I^{-2}(\mu;\lambda)\frac{e^{-\lambda
    x}}{1-e^{-\lambda x}}x^2\,d\lambda
\end{equation}
and a numeric procedure should be used to find the optimal measure for
any given values of $u$ and $n$. Similarly to the case of objective
function $G_3$ above, numeric experiments show that the optimal
measure is one atomic, although we cannot show this rigorously. In
the class of one-atomic measures $\mu=n\delta_{x}$ the goal function
simplifies to
\begin{equation}
\label{eq:G4(x)}
G_4(n\delta_{x};u)=-\frac{e^{ux}-e^{x}-x(u-1)}{nx^3(u-1)}
\end{equation}
Similarly to $G_3$ above, the point of maximum $x_{max}$ of this
function approaches 0 when $u$ grows, see the right plot on
Figure~\ref{fig:G4}. So that when $u$ is such that $x_{max}>1/n$ the
optimal dose is $1/n$ (the whole substrate is used), otherwise the
optimal dose is $x_{max}$. For example, for $n=30$ mice,
the optimal dose is also given by \eqref{eq:G3n30} but with
$u^*\approx 64.47$ this time.

\paragraph{Gamma prior distribution.} Gamma prior distribution arises
naturally as a posterior distribution for the Poisson parameter
$\lambda$ which has a Uniform prior distribution. Since the
probability of having $k$ pre-HSCs in the substrate is
\begin{equation}
\label{eq:Poisson-distribution-function-unconditional-on-lambda}
\P\{\omega(V) =k\}=\E _Q\P\{\omega(V)=k\cond\lambda\}
=\frac{1}{u-1}\int_{1}^{u}\frac{\lambda^{k}e^{-\lambda}}{k!}\,d\lambda,
\end{equation}
The posterior p.d.f.\ for $\lambda$ is then
\begin{equation}
\label{eq:posterior-distribution-for-lambda}
f_{\lambda}(x\cond \omega(V) =k)
=\frac{x^ke^{-x}}{\int_{1}^{u} \lambda^k e^{-\lambda} d\lambda}\propto x^{k}e^{-x},
\end{equation}
which is close to $\Gamma(k+1,1)$ distribution for large $u$. So
once an estimate for the total number of pre-HSCs k is available, the
Gamma prior would be a reasonable for $\lambda$ in the subsequent
experiment. Moreover, posterior distribution for a Gamma prior
$\Gamma(\alpha,1)$ given $k$ pre-HSCs will be $\Gamma(\alpha+k,2)$ and
so on. So generally, $\lambda\sim\Gamma(\alpha,\beta)$ with rate parameter
$\beta\in\N$ could be a reasonable assumption. 

Under this assumption, Criterion~\eqref{eq:G_2} becomes
\begin{equation}
\label{eq:G2(m) with a gamma prior distr}
G_2(\mu ;\alpha ,\beta)=\int_{0}^{\infty}\int_{(0,1]} f_{\alpha,\beta}(\lambda) \frac{e^{-\lambda
    x}}{1-e^{-\lambda x}}x^2 \, \mu(dx) d\lambda\ ,
\end{equation}
where
\begin{equation}
\label{eq:Gamma-density-function}
f_{\alpha,\beta}(\lambda)=\frac{\beta ^{\alpha}}{\Gamma(\alpha)}\lambda^{\alpha -1}e^{-\beta\lambda}.
\end{equation} 
So $G_2$ is a linear functional with the gradient function
\begin{multline}
\label{eq:gradient function of G_2(m) w.r.t m when lambda has gamma distr.}
g_2(x,\mu ; \alpha , \beta)=\int_{0}^{\infty}\frac{e^{-\lambda x}}{1-e^{-\lambda x}}x^2
f_{\alpha,\beta}(\lambda) d\lambda=
\int_{0}^{\infty}r(\lambda x) \lambda^{-2} f_{\alpha,\beta}(\lambda)
d\lambda\\
=\frac{\beta^2}{(\alpha-1)(\alpha-2)} \int_{0}^{\infty}r(\lambda x) f_{\alpha-2,\beta}(\lambda)
d\lambda \\ 
=\frac{\beta^2}{(\alpha-1)(\alpha-2)} \int_{0}^{\infty}r(\lambda) f_{\alpha-2,\beta/x}(\lambda)
d\lambda.
\end{multline}
Thus the situation here is similar to the case of the Uniform
distribution: depending on whether the maximum of this function is
below or above $1/n$ one needs to take equal doses of volumes
$x_{max}$ the point where $g_2$ attains its maximum or $1/n$,
respectively.

\begin{figure}
\includegraphics[width=6cm]{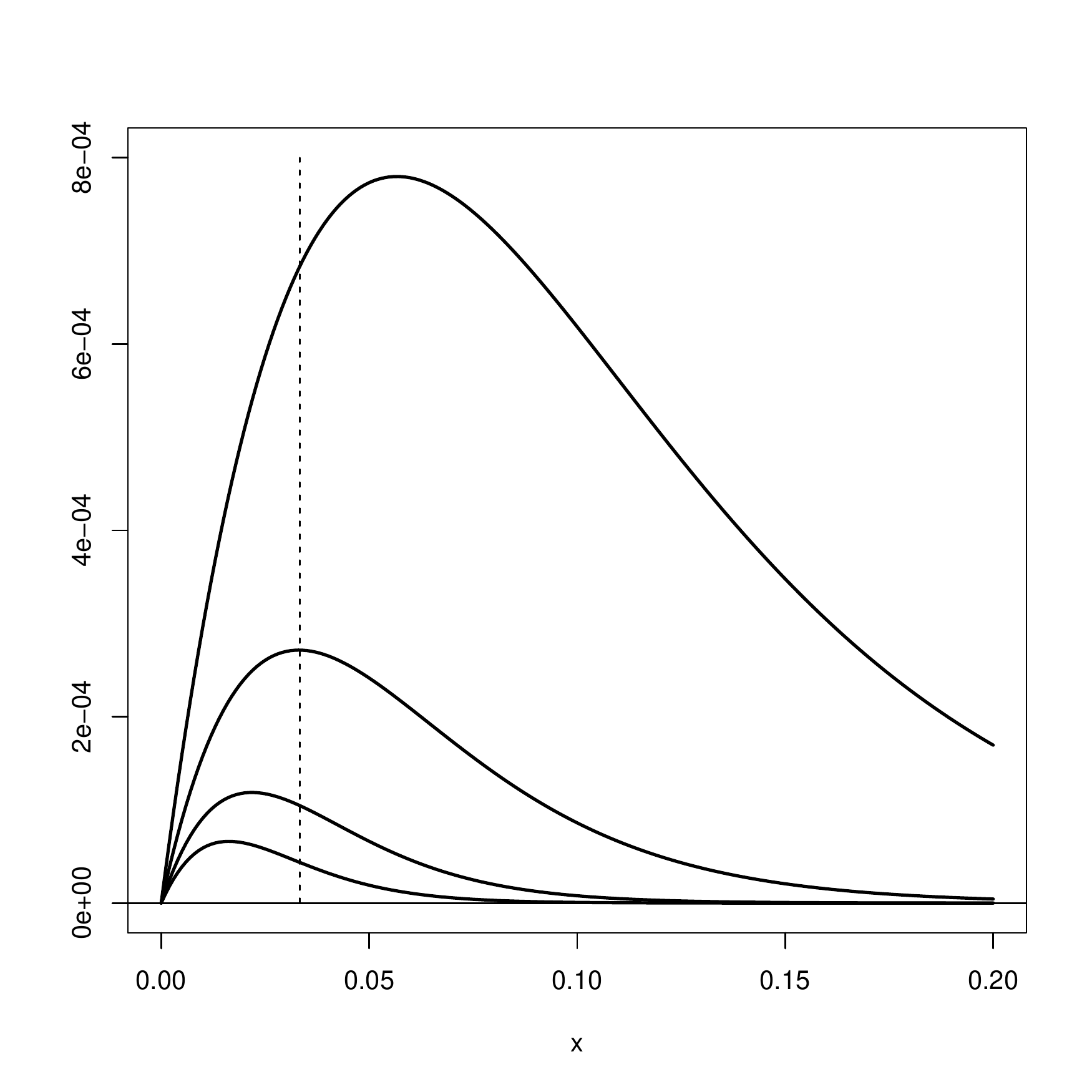}
\includegraphics[width=6cm]{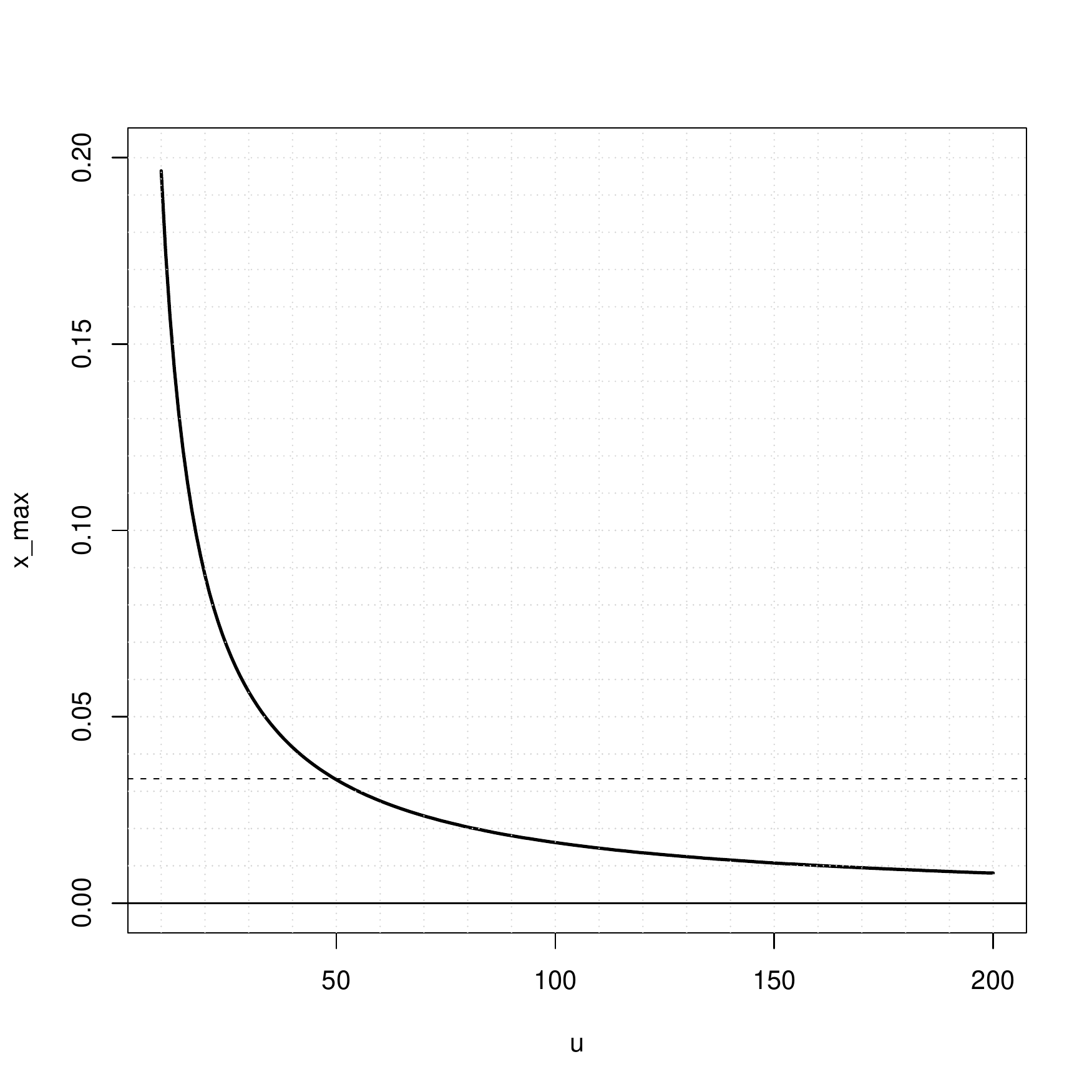}
\caption{Left plot: function $g_2(x)$ for $\beta=1$ and
  $\alpha=30,50,75,100$ (from top to bottom). The vertical dashed line is through the
  point $1/n$ with $n=30$.  Right plot: $x_{max}$ as a function of $\alpha$ (horizontal dashed
  line is through $1/n$ for $n=30$ intersection the curve at point
  $u^*\approx 49.68$).
  \label{fig:g2-gamma}}
\end{figure}

Consider as an example the first-iteration case when $\beta=1$. The
gradient function and the point of maximum are shown in
Figure~\ref{fig:g2-gamma}. The optimal design is given by
\begin{equation}\label{eq:G2n30}
x^*=
 \begin{cases}
 \frac{1}{n} & \quad \alpha\leq \alpha^*;\\
 x_{max}  & \quad \alpha>\alpha^*.  
 \end{cases}
 \end{equation}
where $\alpha^*\approx 49.68$.

\begin{figure}
\includegraphics[width=6cm]{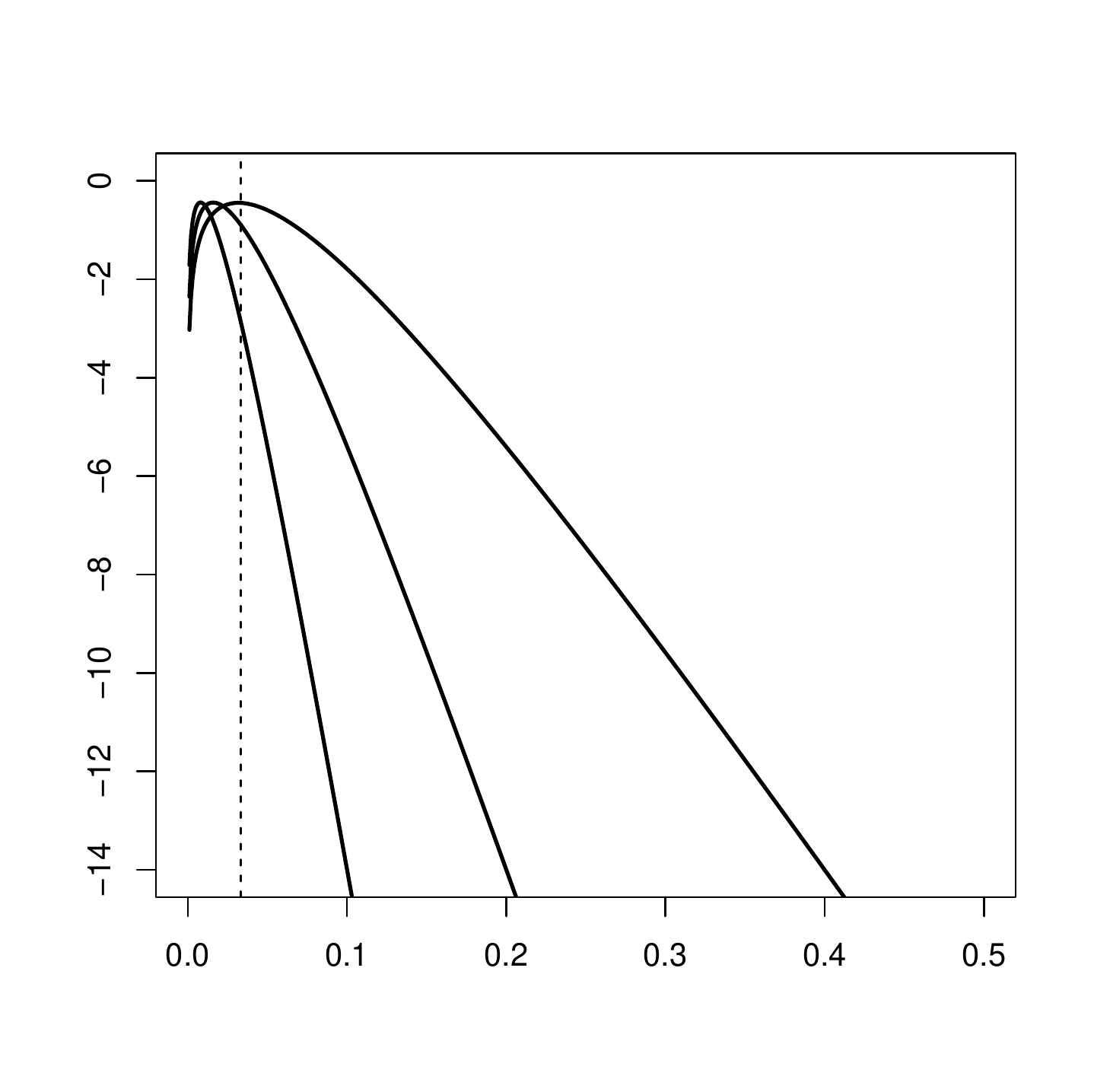}
\includegraphics[width=6cm]{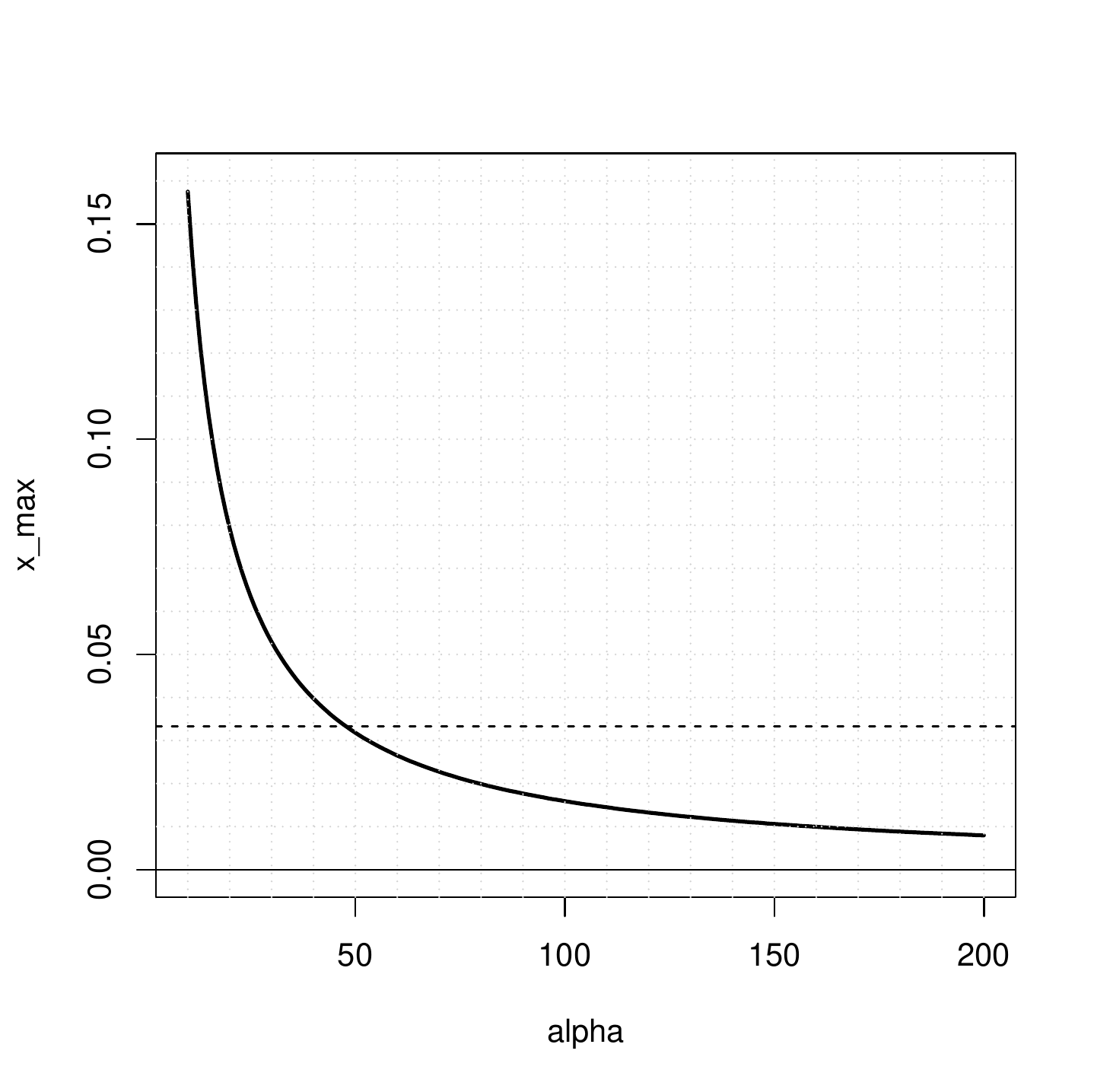}
\caption{Left plot: function ${G}_3(n\delta_x;\alpha)$ for $\beta=1$
  and $\alpha=50, 100$ and 200. Similarly to the Uniform prior case
  as $\alpha$ grows, the point of maximum $x_{max}$ approaches 0. The
  vertical dashed line is through the point $1/n$ with $n=30$.  Right
  plot: $x_{max}$ as a function of $\alpha$ (horizontal dashed line is
  through $1/n$ for $n=30$).\label{fig:G3_Gamma}}
\end{figure}

Under Gamma prior~\eqref{eq:Gamma-density-function}, criterion
function~\eqref{eq:G_3} takes the following form:
\begin{equation}
\label{eq:G_3-Gamma}
G_3(\mu;\alpha,\beta)=\int_{0}^{\infty}f_{\alpha,\beta}(\lambda)\log I(\mu;\lambda) d\lambda ,
\end{equation}
where $I(\mu;\lambda)$ is given by~\eqref{eq:Idef}. Note that $G_3$ is
not a linear functional of a measure, however it is Fr\'echet
differentiable with the gradient function
\begin{equation}
\label{eq:g_3-Gamma}
g_3(x,\mu;\alpha,\beta)=\int_{0}^{\infty}I^{-1}(\mu;\lambda)\frac{e^{-\lambda
    x}}{1-e^{-\lambda x}}x^2 f_{\alpha,\beta}(\lambda) d\lambda. 
\end{equation}
Thus, as before numeric methods for $\beta=1$ and a given values of
$\alpha$ and $n$ should be employed. Our experiments show that the
optimum measure given by the steepest descent algorithm still contains
only one atom. So, optimising~\eqref{eq:G_3-Gamma} for one atomic
measures $n\delta_x$ over $x\in(0,1/n]$ leads to the design given
by~\eqref{eq:G2n30}, with
$\alpha^{*}\approx47.70$. Figure~\ref{fig:G3_Gamma} shows
$G_3(n\delta_{x};\alpha)$ and the optimal equal doses for different values
of $\alpha$.

\begin{figure}
\includegraphics[width=6cm]{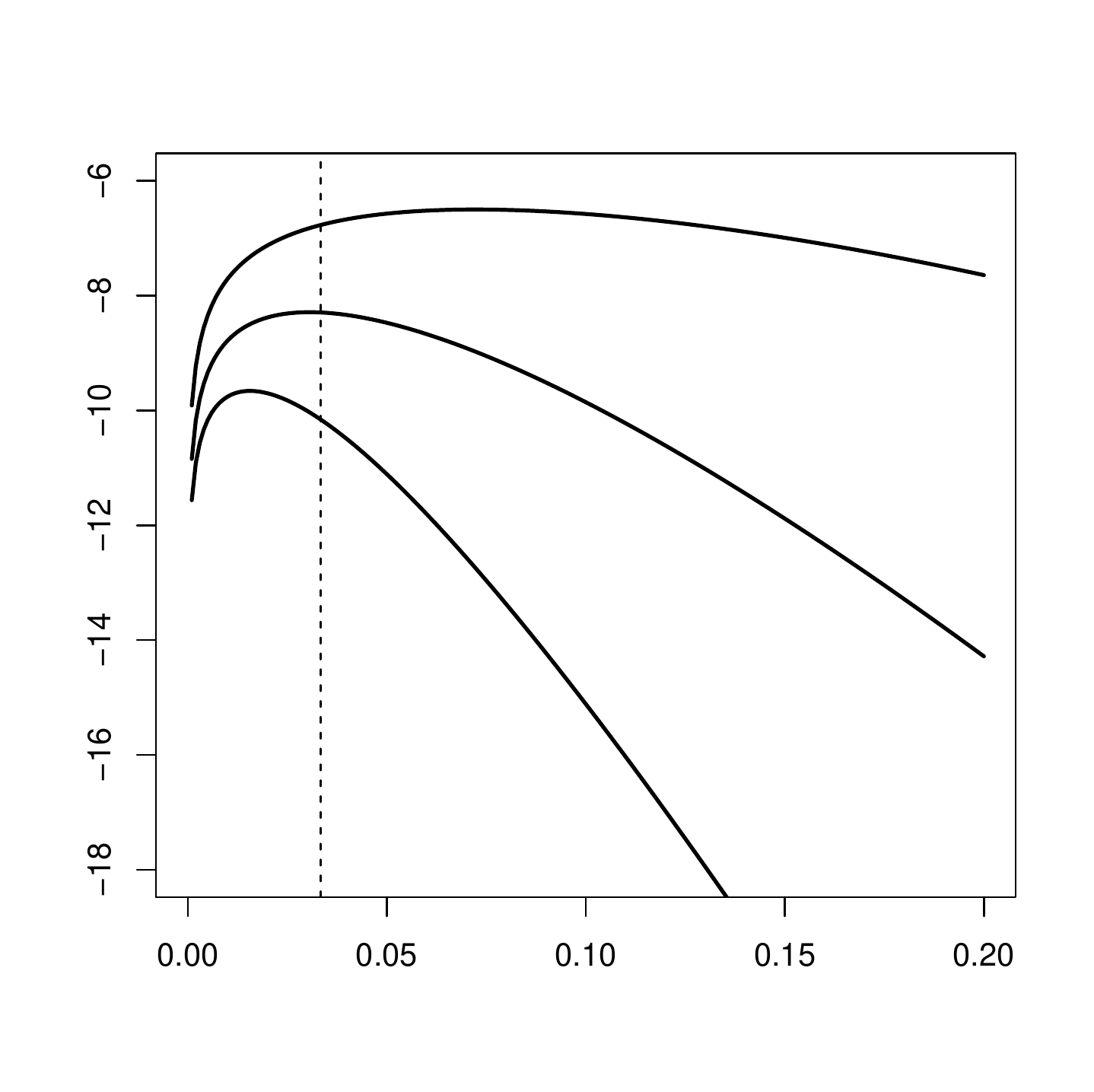}
\includegraphics[width=6cm]{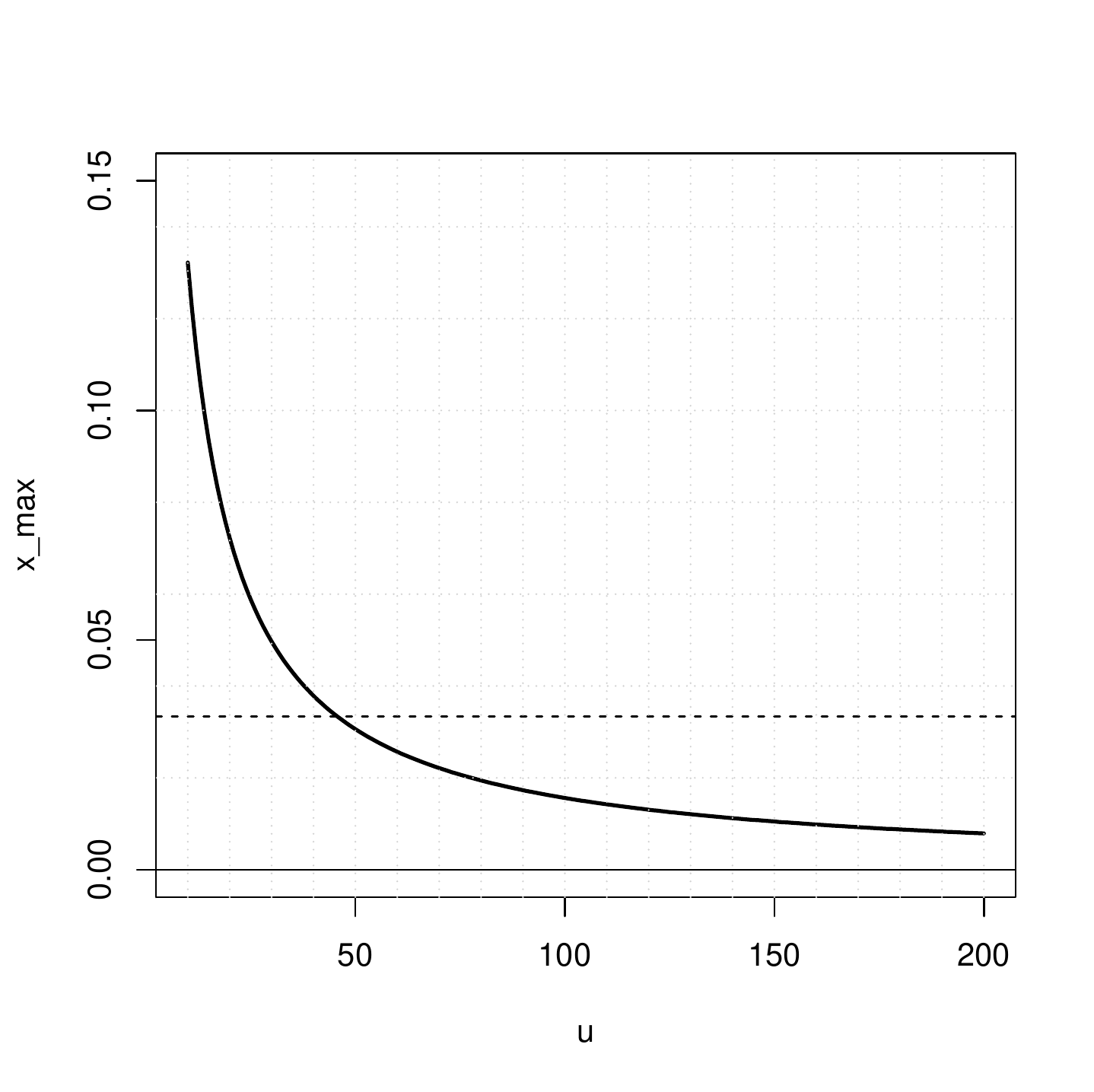}
\caption{Left plot: function $-\log(|nG_4(n\delta_{x};\alpha)|)$ for
  $\beta=1$ and $\alpha=20, 50$ and 100. Similar to the Uniform case
  for $G_4$, as $\alpha$ grows, the point of maximum $x_{max}$
  approaches 0. The vertical dashed line is through the point $1/n$
  with $n=30$.  Right plot: $x_{max}$ as a function of $\alpha$
  (horizontal dashed line is through $1/n$ for
  $n=30$).\label{fig:G4_Gamma}}
\end{figure}

Finally, consider $G4$ in~\eqref{eq:G_4} with Gamma prior distribution on $\lambda$:
\begin{equation}
\label{eq:G_4_Gamma}
G_4(\mu;\alpha,\beta)=-\int_{0}^{\infty}I^{-1}(\mu;\lambda)f_{\alpha,\beta}(\lambda) d\lambda,
\end{equation}
which is Fr\'echet differentiable with the gradient function
\begin{equation}
\label{eq:g_4_Gamma}
g_4(x,\mu;\alpha,\beta)=\int_{0}^{\infty}I^{-2}(\mu;\lambda)\frac{e^{-\lambda
    x}}{1-e^{-\lambda x}}x^2 f_{\alpha,\beta}(\lambda) d\lambda. 
\end{equation}
Here also, our numeric experiments for
$\beta=1$ and any given values of $\alpha$ and $n$, show that the
optimal measure contains a single atom. Thus, in the class of one atomic
measures, \eqref{eq:G_4_Gamma} simplifies to
\begin{equation}
\label{eq:G4_oneatom_Gamma}
G_4(n\delta_{x};\alpha)=-\dfrac{1}{nx^2}\int_{0}^{\infty}\frac{1-e^{-\lambda
    x}}{e^{-\lambda x}} f_{\alpha}(\lambda) d\lambda. 
\end{equation}
which leads to the one-atom optimal design~\eqref{eq:G2n30} with
$\alpha^*\approx 45.74$, see Figure~\ref{fig:G4_Gamma}.

\subsection{Optimal design with additional costs}
\label{subsec:cost}

The measure optimisation formalism we used above allows for seamless
inclusion of different costs associated with
the experiment either as additional terms in the goal function or as
additional contraints. This can be done in both non-Bayesian and
Bayesian settings. Typical loss incured in the above described
experiment are the non-repopulated laboratory mice which are quite
expensive. Added to this, the cost of all the experiment when all mice
die or all mice survive so if it does not yield meaningful
results. In this section we demonstrate how these costs affect the
optimal experimental design.

\paragraph{Non-Bayesian setting.} To take into account the cost
of non-repopulated mice, note that for a fixed $\lambda$, the
indicator $\chi_{i}$ that the mouse which received the $i$th dose
of volume $x_i$ has not repopulated is a Bernoulli random variable
with a parameter $e^{-\lambda x_i}$,
see~\eqref{eq:Bernoulli-random-variable}. 
Thus, the mean number of non-repopulated mice $\E\sum_{i=1}^{n} \chi_{i}$ is
\begin{equation}
\label{eq:T1(x)}
T_1(\mathbf{x} ;\lambda)=\sum_{i=1}^{n} e^{-\lambda x_i}, 
\end{equation} 
which can be represented as a functional of the design measure $\mu$ as follows:
\begin{equation}
\label{eq:T1(mu)}
T_1(\mu ;\lambda)=\int e^{-\lambda x}\mu(dx).
\end{equation}
Recall that as in the previous sections, all the integrals are taken
over the range $(0,1]$. To include the cost associated with a spoilt experiment, compute the
probability that all the mice do not repopulate which can be written as
follows:
\begin{align}
\label{eq:T20(x)}
T_{21}(\mathbf{x};
\nonumber
\lambda)=\prod_{i}^{n} \P\{\chi_{i}=1\}\\
&=\prod_{i=1}^{n}e^{-\lambda x_i}=e^{-\lambda \sum_{i=1}^{n}x_i}.
\end{align}
Consequently, the probability that all the mice repopulate is
\begin{equation}
\label{eq:T21(x)}
T_{20}(\mathbf{x};\lambda)=\prod_{i=1}^{n}\P\{\chi_{i}=0\}
=\prod_{i=1}^{n}(1-e^{-\lambda x_i})=e^{\sum_{i=1}^{n}\log(1-e^{-\lambda x_{i}})}.
\end{equation}
In terms of the design measure, \eqref{eq:T20(x)}
and~\eqref{eq:T21(x)} are the following functionals:
\begin{align}
\label{eq:T20(mu)}
T_{21}(\mu;\lambda) & =e^{-\lambda\int x\mu(dx)}\\
\label{eq:T21(mu)}
T_{20}(\mu;\lambda) & =e^{\int \log(1-e^{-\lambda x})\mu(dx)}.
\end{align}
This gives rise to the expression for the probability of a spoilt
experiment due to either all the mice repopulating or all the mice
non-repopulating:
\begin{equation}
\label{eq:T2(mu)}
T_{2}(\mu;\lambda)=T_{21}(\mu;\lambda)+T_{20}(\mu;\lambda).
\end{equation}
We note that $T_{1}$ and $T_{2}$ are Fr\'echet differentiable with the
gradient functions
\begin{equation}
\label{eq:T1-grad}
t_1(x;\lambda)=e^{-\lambda x}
\end{equation}
and
\begin{equation}
\label{eq:T2-grad}
t_2(x,\mu;\lambda)=
-\lambda x\, T_{21}(\mu ;\lambda)+\log(1-e^{-\lambda x})\, T_{20}(\mu ;\lambda).
\end{equation}

Now consider maximisation of a new goal function
\begin{equation}
\label{eq:G1-tilde}
\tilde{G_1}(\mu;\lambda)=G_1(\mu ;\lambda)-c_1 T_1(\mu ;\lambda)-c_2 T_2(\mu ; \lambda), 
\end{equation}
where $G_1$ is given by~\eqref{eq:G1}, under
constraints~\eqref{eq:tot-mass} and~\eqref{eq:tot-volume}. The positive constants $c_1, c_2$
should be set by the experimentator to reflect the cost of mice and of
a spoilt experiment which should be offset against the usefulness of
the results reflected in the original goal function $G_1$ involving
the Fisher information.

The function $\tilde{G_1}$ is strongly diferentiable 
with the gradient function given by
\begin{equation}
\label{eq:G1-tilde-grad}
\tilde{g_1}(x,\mu;\lambda)=g_1(x;\lambda)-c_1t_1(x;\lambda)-c_2t_2(x,\mu;\lambda),
\end{equation}
where $g_1$ as in~\eqref{eq:gradient-of-G_1(m)}.

\begin{figure}[ht]
  \centering
  \includegraphics[width=7cm]{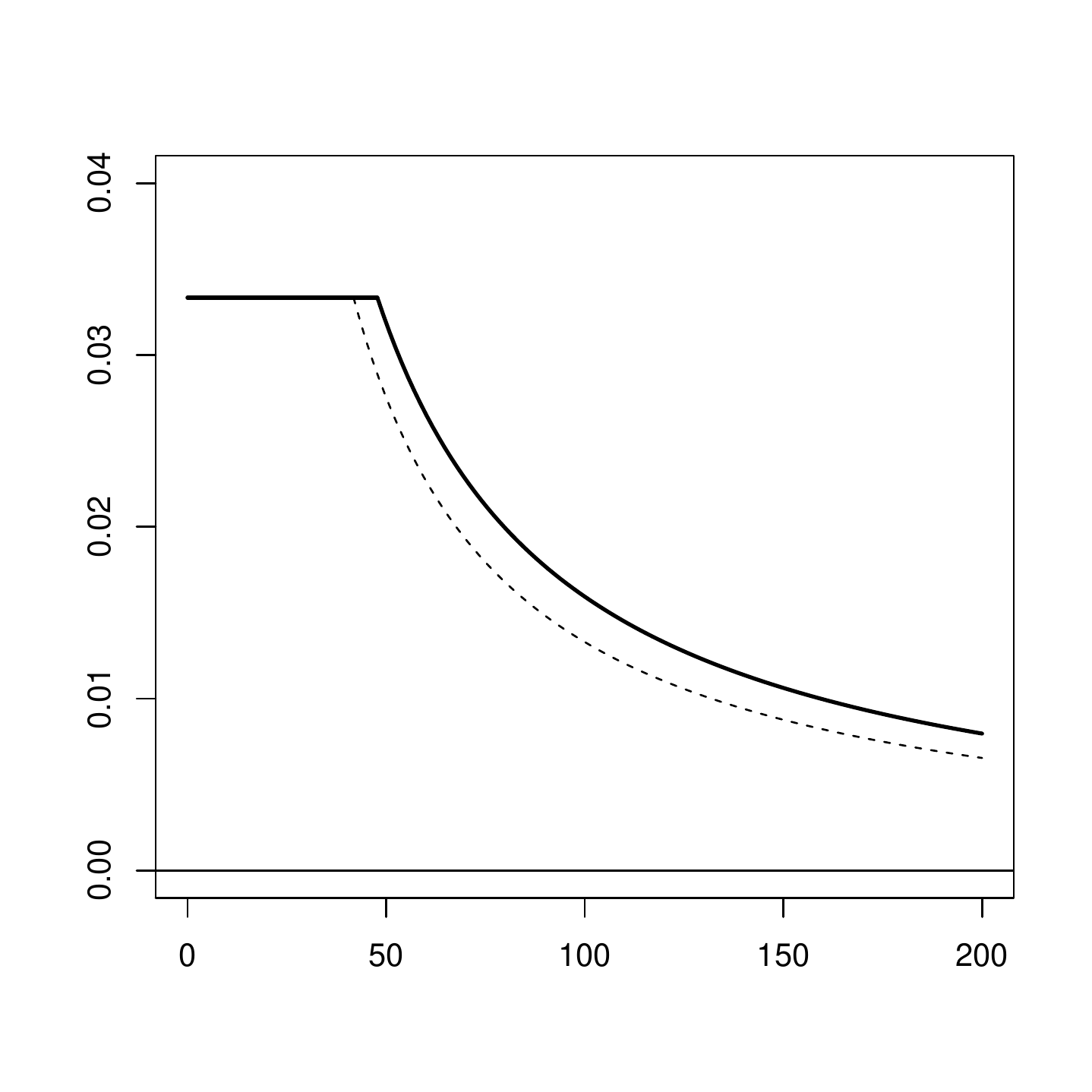}
  \caption{ The optimal doses for different values of $\lambda$ for
    the goal function $\tilde{G_1}$ when $n=30$, $c_1=10^{-4}$ and
    $c_2=1$ (dashed) and for $G_1$ (solid line) as on Plot~(c) in
    Figure~\ref{fig:r}.}
  \label{fig:tildeG1-dose-lam}
\end{figure}

Figure~\ref{fig:tildeG1-dose-lam} shows the optimal dose as a function
of $\lambda$, for $n=30$, $c_1=10^{-4}$ and $c_2=1$.  Compared to the
optimal design for $G_1$ shown on Figure~\ref{fig:r} one has to dilute
more and also to start dilution earlier at $\lambda=41.8$ compared to
47.8 for $G_1$ which reflects the caution not to allow all the mice to
survive as this would mean a spoilt experiment.

\paragraph{Bayesian setting.} Under $\Gamma(\alpha,1)$ prior
distribution for $\lambda$, the average number of non-repopulated mice
becomes
\begin{equation}
\label{eq:T1-Gamma-Prior}
\E_{Q} T_1(\mu ;\lambda)=\int\frac{1}{(1+x)^{\alpha}}\mu(dx),
\end{equation}
with the gradient function 
\begin{equation}
\label{eq:T1-grad-Gamma-Prior}
\E_{Q}t_1(x; \lambda)=\frac{1}{(1+x)^{\alpha}}.
\end{equation}
The probability of a spoilt experiment is
\begin{align}
\label{eq:T2-Gamma-Prior}
\E_{Q}T_2(\mu ;\lambda)&=\E_{Q}T_{21}(\mu ;\lambda)+\E_{Q}T_{20}(\mu ;\lambda)\\
&=\dfrac{1}{(1+H(\mu))^{\alpha}}+\E_{Q}T_{20}(\mu ;\lambda),
\end{align}
with the corresponding gradient function
\begin{equation}
\label{eq:T2-grad-Gamma-Prior}
\E_{Q}t_2(x,\mu;\lambda)=-\dfrac{\alpha x}{(1+H(\mu))^{\alpha +1}} 
+\E_{Q}\log(1-e^{-\lambda x})T_{20}(\mu ;\lambda),
\end{equation}
where $H(\mu)=\int x\mu(dx)$.

Take $G_4$ as in~\eqref{eq:G_4_Gamma} with $\beta=1$ and consider a
new goal function
\begin{equation}
\label{eq:G4-tilde}
\tilde{G_4}(\mu ;\alpha)=G_4(\mu ;\alpha)-c_1 \E_{Q}T_1(\mu ;\lambda)-c_2\, \E_{Q} T_2 (\mu ;\lambda).
\end{equation}
It is Fr\'echet differentiable and possesses a gradient function
\begin{equation}
\label{eq:G4-tilde-grad}
\tilde{g_4}(x,\mu;\alpha)=g_4(\mu ,x ;\alpha)-c_1\E_{Q}t_1(x;\lambda)-c_2\E_{Q}t_2(x,\mu;\lambda).
\end{equation} 
We are aiming to maximise $\tilde{G_{4}}$ under the same constraints
\eqref{eq:tot-mass} and~\eqref{eq:tot-volume}.

\begin{figure}[ht] \centering
  \includegraphics[width=7cm]{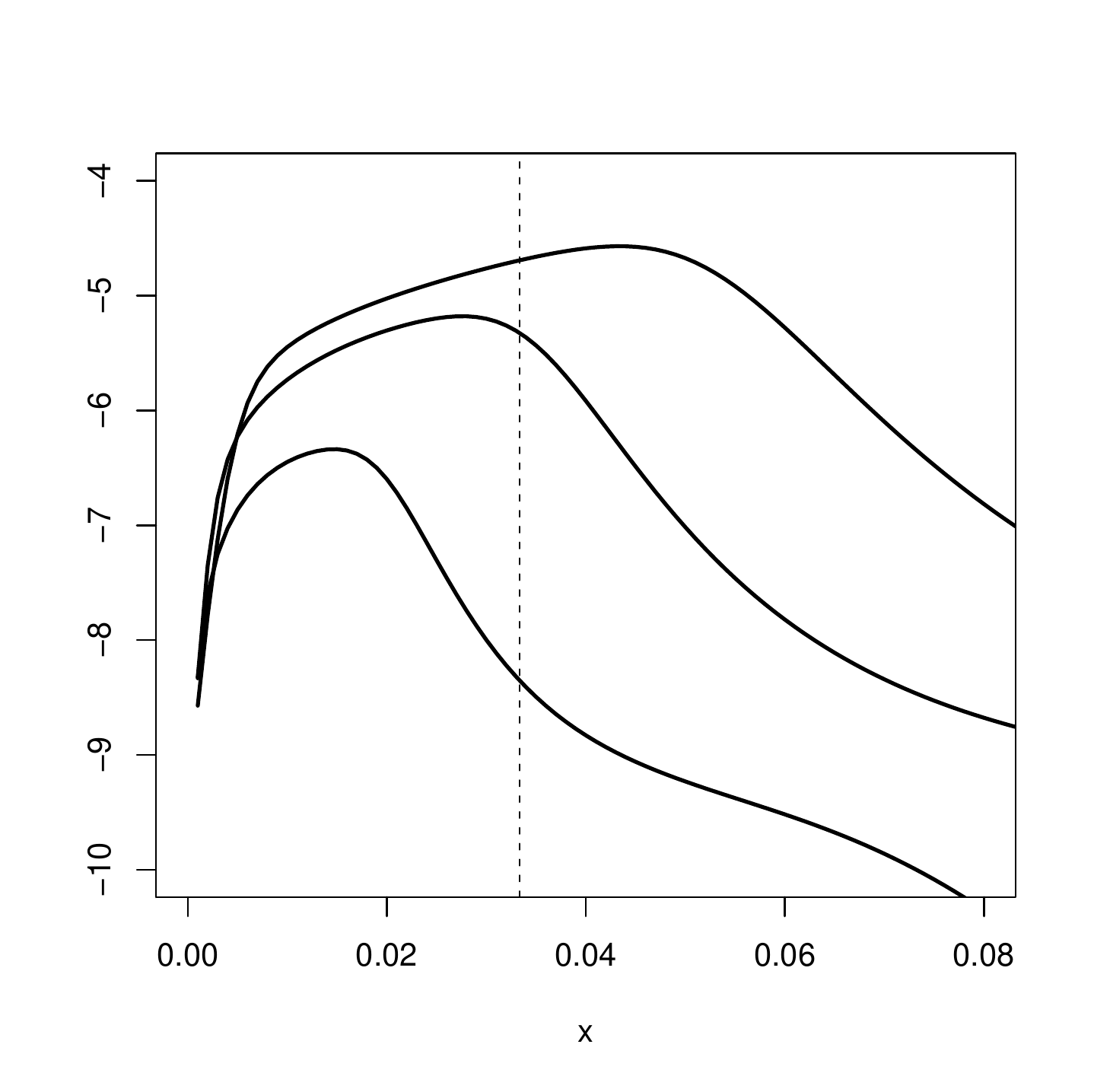}
  \caption{Function $-\log|\tilde{G_4}(n\delta_{x};\alpha)|$ for the
    costs $c_1=0.005$, $c_2=5$ and values $\alpha=30,50,100$ (in the
    order of a decrease). The vertical dashed line is through the
    point $1/n$ with $n=30$. The extremal points are $x_{max}=0.043,
    0.028$ and $0.015$ leading to optimal designs $1/30, 0.028$ and
    $0.015$, respectively.}
  \label{fig:GG4-one-atom-Gamma}
\end{figure}

Again, numeric solutions are still one-atomic for different tried
values of the parameters $\alpha, c_1$ and $c_2$ we tried.
Figure~\ref{fig:GG4-one-atom-Gamma} shows a typical picture for
$-\log|\tilde{G_4}(n\delta_{x};\alpha)|$ for the class of one-atomic
design measures and the optimal doses for three different values of
$\alpha$.

\paragraph{Mixture prior distribution}
Finally, consider the setting typical for testing two alternative
hypotheses about possible values of the mean number of
pre-HSCs. Suppose that $\lambda$ can take two different values, our prior belief is that it is rather
$\lambda_1$ than $\lambda_2$ with probability $0<p<1$. This gives rise
to optimisation of the following goal function
\begin{equation}
\label{eq:G1M}
G_{1M}(\mu ;p,\lambda_{1},\lambda_{2})=p\ G_{1}(\mu;\lambda_{1})+(1-p)\ G_{1}(\mu;\lambda_{2}).
\end{equation}
under constraints \eqref{eq:tot-mass} and \eqref{eq:tot-volume}.

\begin{figure}[ht]
  \centering
  \includegraphics[width=12cm]{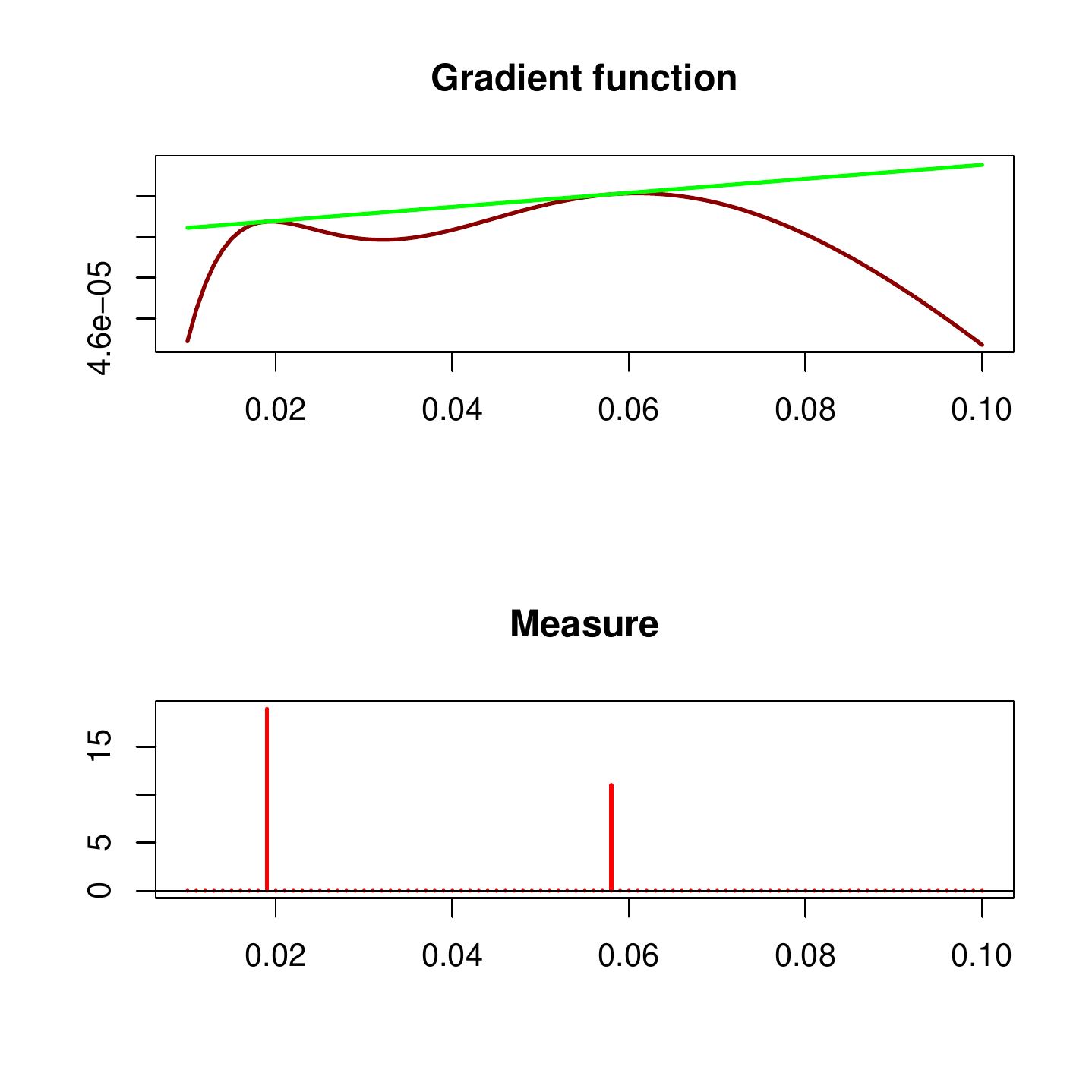}
  \caption{Solution obtained by \texttt{medea} for $p=0.05$,
    $\lambda_1=25$, and $\lambda_2=150$ shows that the numeric
    solution is a two atomic measure $18.98\ \delta_{0.019}+11.02\
    \delta_{0.058}$ and all the substrate should be used.}
  \label{fig:g1M-medea.pdf}
\end{figure}

The gradient function for~\eqref{eq:G1M} is 
\begin{equation}
\label{eq:G1M-grad}
g_{1M}(x;p,\lambda_{1},\lambda_{2})=p\ g_{1}(x;\lambda_{1})+(1-p)\ g_{1}(x;\lambda_{2}).
\end{equation} 
and numerical optimisation can be employed to find the optimal design
for any particular values of $\lambda_1, \lambda_2$ and $p$.  For
example of $\lambda_1=25, \lambda_2=150$ and $p=0.05$,
Figure~\ref{fig:g1M-medea.pdf} shows a numerically obtained optimal
measure which is \emph{two-atomic}: 19 mice should receive equal
doses of $0.019$ and the rest 11 of them should receive the same doses of
$0.058$.

\section{Discussion}
\label{sec:discussion}
We have considered dilution experiments with volume constraints
typical in biological and medical research. An important particularity
of the experiments we consider here is that the time delay necessary
for pre-HSCs to develop and then for the injected doses to take effect
in mice prevents from planing multiple stage experiments: the
estimation should necessarily be done from the first and only stage.
As we have seen, in all the considered cases of the goal functions, the
optimal design is attained on a one-point measure, meaning that all
the doses should have equal volume. This parallels the well known
result about the D-optimal design measure for a linear regression
model: The Kiefer-Wolfowitz theorem assures that such measure is
atomic with the number of atoms to be at most the number of contraints
plus one, \seg\cite{Kiefer:1960}. In our case we have one
volume constraint (in addition to the measure to have a fixed total
mass) so the number of atoms is at most two. Indeed, the gradient
functions we observed are convex on the interval from 0 to the point of
maximum, so the only possibility to satisfy the necessary optimality
condition given in Theorem~\ref{th:optimisation-with-a-limited-cost},
is for the optimal measure to have only one atom. However, the example
of a mixture prior distribution considered in the last section shows
that two-atom designs are indeed possible.

In practical terms, if the number of pre-HSC expected to
be not very large (a few dozens or less), the whole substrate should
be used to derive the doses, otherwise only part of it. We have
characterised above what is `not very large' and how much the doses
should be diluted. We have fixed some of the parameters here (the
number of mice $n=30$, $\beta=1$), driven by
practical applications to stem cell research. For other applications
and other values of these parameters one make use the computer codes
the authors make freely available for download.

Advantage of measure optimisation approach is that additional
requirements can be easily incorporated into the goal function or
added as further constraints, \eg limited cost associated with
non-repopulated mice and/or the cost of the whole experiment if all
mice are repopulated or all are not repopulated.


\end{document}